\documentclass[twoside,11pt,abbrvbib]{article}

\usepackage{jmlr2e}

\usepackage{natbib}
\usepackage[utf8]{inputenc} 
\usepackage[T1]{fontenc}    
\usepackage{lmodern}
\usepackage{booktabs}       
\usepackage{amsfonts}       
\usepackage{nicefrac}       
\usepackage{microtype}      

\usepackage{amsmath}
\usepackage{graphicx}
\usepackage{epstopdf}
\usepackage{subcaption}
\usepackage{multirow}
\usepackage{booktabs}
\usepackage{lineno}
\usepackage{float}
\usepackage{setspace}
\usepackage{overpic}
\usepackage{rotating}
\usepackage{wrapfig,lipsum}
\usepackage{mathtools}
\usepackage{bm}

\usepackage[format=plain,font=small,margin=0.0in,labelfont={sc},labelsep={colon}]{caption}

\usepackage{bibspacing}
\usepackage{tabulary}

\usepackage{algorithm}
\usepackage{algpseudocode}
\algnewcommand\algorithmicinput{\textbf{Input:}}
\algnewcommand\Input{\item[\algorithmicinput]}

\usepackage{tikz}
\usetikzlibrary{positioning,shapes,shadows,arrows,calc}

\tikzstyle{block} = [draw,line width=1.2pt, fill=white!30, rectangle, 
minimum height=3em, minimum width=3em, rounded corners]         
\tikzstyle{blocky} = [draw, line width=1.2pt,fill=yellow!40, rectangle, 
minimum height=2.5em, minimum width=2.5em, rounded corners]    
\tikzstyle{blockp} = [draw,line width=1.2pt, fill=purple!30, rectangle, 
minimum height=3em, minimum width=3em, rounded corners]       
\tikzstyle{blockb} = [draw,line width=1.2pt, fill=blue!20, rectangle, 
minimum height=2.5em, minimum width=2.5em, rounded corners]     
\tikzstyle{blockr} = [draw, line width=1.2pt,fill=red!30, rectangle, 
minimum height=2.5em, minimum width=2.5em, rounded corners]                 
\tikzstyle{blockgr} = [draw, line width=1.2pt,fill=black!10, rectangle, 
minimum height=2.5em, minimum width=2.5em, rounded corners]                  
\tikzstyle{line}=[-, line width=1.2pt]

\tikzstyle{blocky4} = [draw, line width=1.2pt,fill=white!10, rectangle, 
text centered,minimum height=5.0em, text width=6.5em, rounded corners]        

\tikzstyle{blockyY} = [draw, line width=1.2pt,fill=red!10, rectangle, 
text centered,minimum height=2.5em, text width=6.5em, rounded corners]    

\tikzstyle{blockyZ} = [draw, line width=1.2pt,fill=red!10, rectangle, 
text centered,minimum height=2.5em, text width=4.0em, rounded corners] 

\tikzstyle{blockyC} = [draw, line width=1.2pt,fill=red!10, rectangle, 
text centered,minimum height=2.5em, text width=4.0em, rounded corners] 

\DeclareMathOperator*{\argmin}{arg\rm{}min}

\definecolor{myRed}{RGB}{222,45,38}
\definecolor{darkred}{RGB}{44,162,95}

\newcommand{\bA}{\mathbf{A}}
\newcommand{\bB}{\mathbf{B}}
\newcommand{\bC}{\mathbf{C}}
\newcommand{\bQ}{\mathbf{Q}}
\newcommand{\ba}{\mathbf{a}}
\newcommand{\bb}{\mathbf{b}}
\newcommand{\bc}{\mathbf{c}}

\DeclareMathOperator{\E}{\mathbb{E}}
\newcommand{\tX}{\bm{\mathcal{X}}}
\newcommand{\tY}{\bm{\mathcal{Y}}}

\newcommand{\tB}{\bm{\mathcal{B}}}
\newcommand{\tM}{\bm{\mathcal{M}}}

\newcommand{\tE}{\bm{\mathcal{E}}}

\newtheorem{Theorem}{Theorem}
\newtheorem{Definition}{Definition}

\ShortHeadings{Randomized CP Tensor Decomposition}{Erichson, Manohar, Brunton, Kutz}
\firstpageno{1}

\begin{document}

\title{Randomized CP Tensor Decomposition}

\author{\name N. Benjamin Erichson \email erichson@berkeley.edu \\
       \addr Department of Statistics\\
       University of California, Berkeley, USA
       \AND
       \name Krithika Manohar	\email kmanohar@caltech.edu \\
       \addr Department of Applied \& Computational Mathematics\\
       California Institute of Technology, Pasadena, USA 
       \AND
       \name Steven L. Brunton \email sbrunton@uw.edu  \\
       \addr Department of Mechanical Engineering\\
       University of Washington, Seattle, USA
       \AND       
       \name J. Nathan Kutz \email kutz@uw.edu \\
       \addr Department of Applied Mathematics\\
       University of Washington, Seattle, USA
    }
\editor{}

\maketitle
\sloppy

\begin{abstract}
The \textit{CANDECOMP/PARAFAC} (CP) tensor decomposition is a popular dimensionality-reduction method for multiway data. Dimensionality reduction is often sought after since many high-dimensional tensors have low intrinsic rank relative to the dimension of the ambient measurement space. However, the emergence of `big data' poses significant computational challenges for computing this fundamental tensor decomposition. By leveraging modern randomized algorithms, we demonstrate that coherent structures can be learned from a smaller representation of the tensor in a fraction of the time. Thus, this simple but powerful algorithm enables one to compute the approximate CP decomposition even for massive tensors. The approximation error can thereby be controlled via oversampling and the computation of power iterations. In addition to theoretical results, several empirical results demonstrate the performance of the proposed algorithm. 
\end{abstract}


\section{Introduction}

Advances in data acquisition and storage technology have enabled the acquisition of massive amounts of data in a wide range of emerging applications. In particular, numerous applications across the physical, biological, social and engineering sciences generate large multidimensional,  multi-relational and/or multi-modal data. 
%
Efficient analysis of this data requires dimensionality reduction techniques.
However, traditionally employed matrix decompositions techniques such as the singular value decomposition (SVD) and principal component analysis (PCA) can become inadequate when dealing with multidimensional data. 
This is because reshaping multi-modal data into matrices, or {\em data flattening}, can fail to reveal important structures in the data. 

Tensor decompositions overcome the information loss from flattening. The canonical CP tensor decomposition expresses an $N$-way tensor as a sum of rank-one tensors to extract multi-modal structure. It is particularly suitable for data-driven discovery, as shown by~\citep{hong2020siam} for various learning tasks on real world data.  
However, tensor decompositions of massive multidimensional data pose a tremendous computational challenge. Hence, innovations that reduce the computational demands have become increasingly relevant in this field. 
%
%
The idea of tensor compression, for instance, eases computational bottlenecks by constructing smaller (compressed) tensors, which are then used as a proxy to efficiently compute CP decompositions. Compressed tensors may be obtained, for example, using the Tucker decomposition of a tensor into one small core tensor and $N$ unitary matrices However, this approach requires the expensive computation of the left singular vectors for each mode. 

\paragraph{Related Work.}
This computational challenge can be tackled using modern randomized techniques developed to compute the SVD. 
\cite{tsourakakis2010mach} presents a randomized Tucker decomposition algorithm based on the random sparsification idea of~\cite{achlioptas2007fast} for computing the SVD.  
\cite{zhou2014decomposition} proposed an accelerated randomized CP algorithm using the ideas of~\cite{halko2011rand} without the use of power iterations. An alternative randomized tensor algorithm proposed by~\cite{DrineasTensor} uses random column selection. Related work by \cite{IncompleteTensors} proposed a sparsity-promoting algorithm for incomplete tensors using compressed sensing.

Another efficient approach to building large-scale tensor decompositions is the subdivision of a tensor into a set of blocks. These smaller blocks can then be used to approximate the CP decomposition of the full tensor in a parallelized or distributed fashion~\citep{phan2011parafac}. \cite{ParallelRandomly} fused random projection and blocking into a highly computationally efficient algorithm. More recently, \cite{RandomizedLathauwer} proposed a block sampling CP decomposition method for the analysis of large-scale tensors using randomization, demonstrating significant computational savings while attaining near optimal accuracy. These block based algorithms are particularly relevant if the tensor is too large to fit into fast memory.

\paragraph{Contributions.} 
We present the randomized CP algorithm, which is closely related to the randomized methods of \cite{halko2011rand}.
%
%
Our method proceeds in two stages. The first stage uses random projections with power iterations to obtain a compressed tensor (Algorithm~\ref{alg:qb}). The second stage applies either alternating least squares (ALS) or block coordinate descent (BCD) to the compressed tensor (Algorithms~\ref{alg:als} and~\ref{alg:bcd}), at significantly reduced cost. Finally, the CP decomposition of the original tensor is obtained by projecting the compressed factor matrices back using the basis derived in Algorithm~\ref{alg:qb}.

Randomized algorithms for CP decomposition have been proposed~\citep{zhou2014decomposition,battaglino2018practical}, albeit without incorporating power iterations. 
The power iteration concept is fundamental for modern high-performance randomized matrix decompositions, but to the best of our knowledge, has not been applied in the context of tensors. Without power iterations, the performance of randomized algorithms based on random projections can suffer considerably in the presence of white noise. We combine power iterations with oversampling to further control the error of the decomposition.

\begin{figure}[!t]
	\centering
	\begin{subfigure}[t]{0.92\textwidth}
		\centering
		\DeclareGraphicsExtensions{.pdf}
		\includegraphics[width=1\textwidth]{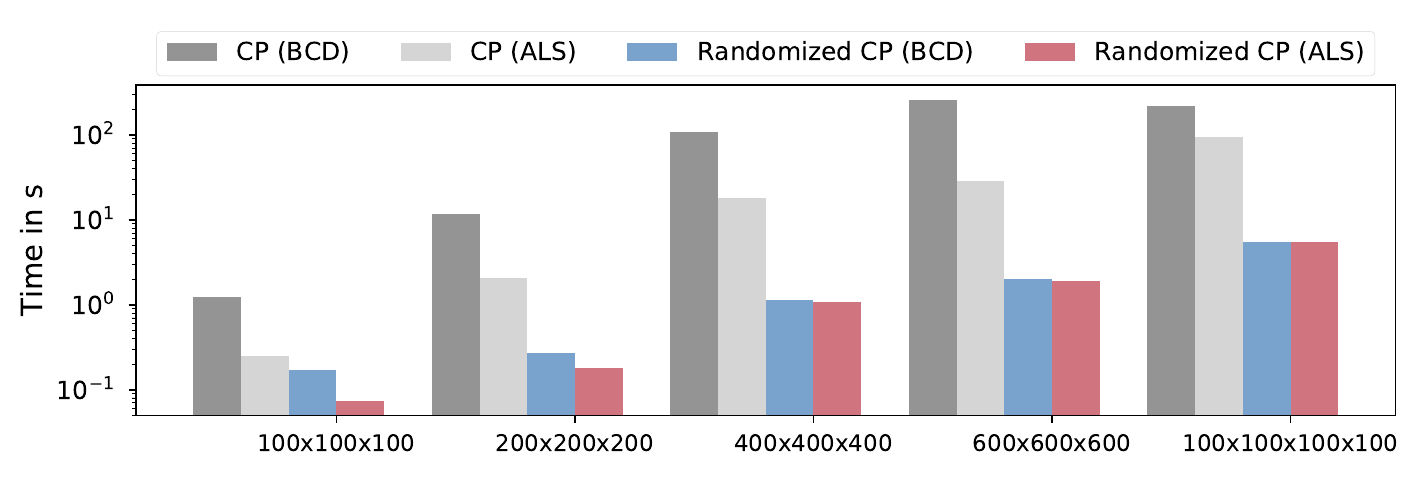}
	\end{subfigure}
	
	\begin{subfigure}[t]{0.92\textwidth}
		\vspace{-0.25in}
		\centering
		\DeclareGraphicsExtensions{.pdf}
		\includegraphics[width=1\textwidth]{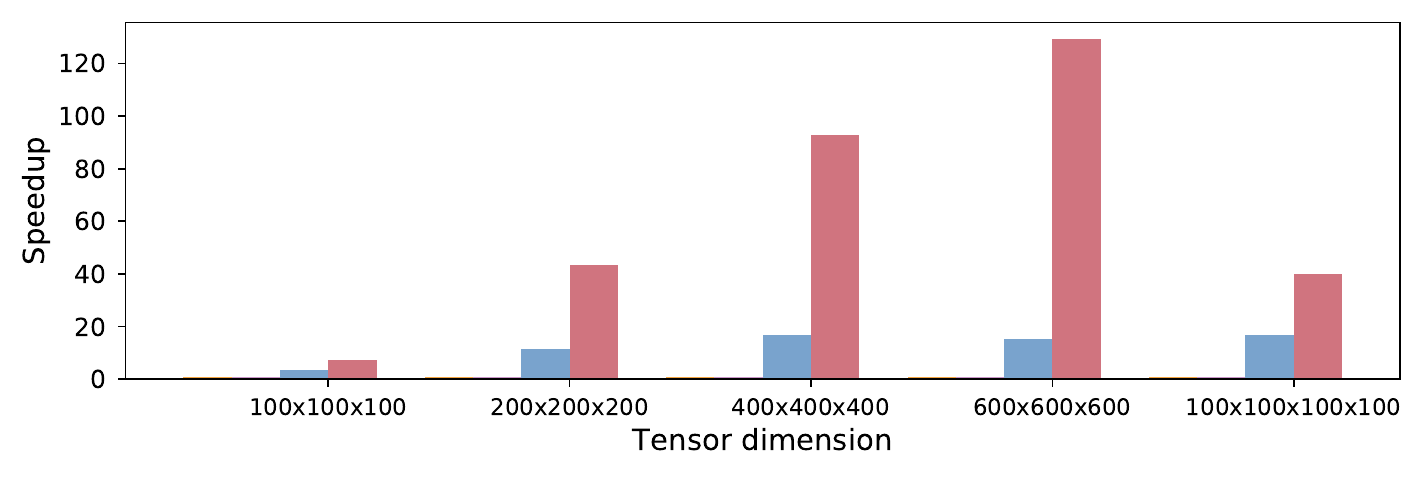}
	\end{subfigure}	
	\vspace*{-.1in}
	\caption{Algorithm runtimes and speedups for varying tensor dimensions for a rank $k=20$ approximation. Speedups rise sharply with increasing dimensions.}
	\label{fig:timing_bar}
\end{figure}

Embedding the CP decomposition into this probabilistic framework allows us to achieve significant computational savings, while producing near-optimal approximation quality.
For motivation, Figure~\ref{fig:timing_bar} shows the computational savings of a rank $k=20$ approximation for varying tensor dimensions. Here we compare the ALS and BCD solver for computing the deterministic and randomized CP decomposition. The computational advantage of the randomized algorithms is significant, while sacrificing a negligible amount of accuracy.      

\paragraph{Organization.} 
The paper is organized as follows: Section~\ref{sec:background} briefly reviews the CP decomposition and randomized matrix decompositions. Section~\ref{sec:cCP} introduces the randomized CP tensor decomposition algorithm.  Section~\ref{sec:results} presents the evaluation of the computational performance, and examples. Finally, Section~\ref{sec:conclusion} summarizes the research findings.

\section{Background} \label{sec:background}
Ideas for multi-way factor analysis emerged in the 1920s with the formulation of the polyadic decomposition by~\cite{hitchcock1927polyadic}.
However, the polyadic tensor decomposition only achieved popularity much later in the 1970s with the canonical decomposition (\textit{CANDECOMP}) in psychometrics, proposed by~\cite{carroll1970candecomp}.  Concurrently, the method of parallel factors (\textit{PARAFAC}) was introduced in chemometrics by~\cite{harshman1970parafac}. Hence, this method became known as the CP (\textit{CANDECOMP/PARAFAC}) decomposition. 
However, in the past, computation was severely inhibited by available computational power. Today, tensor decompositions enjoy increasing popularity, yet runtime bottlenecks persist. 

\subsection{Notation}
Scalars are denoted by lower case letters $x$, vectors as bold lower case letters $\mathbf{x}$, and matrices by bold capitals $\mathbf{X}$. 
Tensors are denoted by calligraphic letters $\tX$. The mode-$n$ unfolding of a tensor is expressed as $\tX_{(n)}$, while the mode-$n$ folding of a matrix is defined as $\mathbf{X}_{(n)}$.
The vector outer product, the Kronecker product, the Khatri-Rao product and the Hadamard product are denoted by $\circ$, $\otimes$, $\odot$, and $\ast$, respectively. 
The inner product of two tensors is expressed as $\langle \cdot ,\cdot \rangle$, and $ \| \cdot \|_F$ denotes the Frobenius norm for both matrices and tensors.

\subsection{CP Decomposition}\label{sec:CP}

\begin{figure}[!b]
	\centering
	\DeclareGraphicsExtensions{.pdf}
	\includegraphics[width=0.9\textwidth]{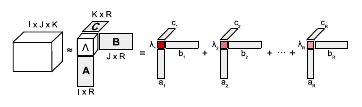}
	\vspace{-0.18in}
	\caption{Schematic of the CP decomposition.}
	\label{fig:CP}
\end{figure}

The CP decomposition is the tensor equivalent of the SVD since it approximates a tensor by a sum of rank-one tensors. Here, {\em tensor rank} is defined as the smallest sum of rank-one tensors required to generate the tensor~\citep{hitchcock1927polyadic}. The CP decomposition approximates these rank-one tensors. 
Given a third order tensor $\tX \in \mathbb{R}^{I \times J \times K}$, the rank-$R$ CP decomposition is expressed as 
\begin{equation}
\tX  \approx \sum_{r=1}^{R} \mathbf{a}_r  \circ  \mathbf{b}_r  \circ  \mathbf{c}_r,
\label{eq:CP}
\end{equation}
where $\circ$ denotes the outer product. Specifically, each rank-one tensor is formulated as the outer product of the rank-one components $\mathbf{a}_r \in \mathbb{R}^{I}$, $\mathbf{b}_r \in \mathbb{R}^{J}$, and $\mathbf{c}_r \in \mathbb{R}^{K}$. Components are often constrained to unit length with the weights absorbed into the vector $\boldsymbol{\lambda}=[\lambda_1, ..., \lambda_R] \in \mathbb{R}^{R}$. Equation~\eqref{eq:CP} can then be re-expressed as (See Fig.~\ref{fig:CP})
\begin{equation}
\tX  \approx \sum_{r=1}^{R}  \lambda_r \cdot \mathbf{a}_r  \circ  \mathbf{b}_r  \circ  \mathbf{c}_r.
\label{eq:CPlamb2}
\end{equation}
More compactly the components can be expressed as factor matrices, i.e., $\mathbf{A}=[\mathbf{a}_1, \mathbf{a}_2,...,\mathbf{a}_R]$, $\mathbf{B}=[\mathbf{b}_1, \mathbf{b}_2,...,\mathbf{b}_R]$, and $\mathbf{C}=[\mathbf{c}_1, \mathbf{c}_2,...,\mathbf{c}_R]$ . Using the Kruskal operator as defined by \cite{kolda2009tensor}, \eqref{eq:CPlamb2} can be more compactly expressed as
\begin{equation*}
\tX  \approx [\![\bm{\lambda; \bA,\bB,\bC}  ]\!].
\label{eq:Tuckeroperator}
\end{equation*}

\subsection{Randomized Matrix Algorithms}\label{sec:background_rand}

The efficient computation of low rank matrix approximations is a ubiquitous problem in machine learning and data mining. Randomized matrix algorithms have been demonstrated to be highly competitive and robust when compared to traditional deterministic methods.
Randomized algorithms aim to construct a smaller matrix (henceforth called sketch) designed to capture the essential information of the source matrix~\citep{Mahoney2011,RandNLA}. The sketch can then be used for various learning tasks. There exist several strategies for obtaining such a sketch, with random projections being the most robust off-the-shelf approach. 
Randomized algorithms have been in particular studied for computing the near-optimal low-rank SVD~\citep{frieze2004fast,liberty2007randomized,woolfe2008fast,Martinsson201147}. 
Following the seminal work by~\cite{halko2011rand}, a randomized algorithm computes the low-rank matrix approximation
\begin{equation*}\label{eq:QB}
\begin{array}{cccc}
\mathbf{A} & \approx & \mathbf{Q} & \mathbf{B} \\[-.05in]
m\times n &   &  m\times k & k\times n
\end{array} 
\end{equation*}
where the target rank is denoted as $k$, and assumed to be $k\ll \text{min}\{m,n\}$.
Intuitively, we construct a sketch $\mathbf{Y} \in \mathbb{R}^{m\times k}$ by forming a random linear weighted combination of the columns of the source matrix $\mathbf{A}$. More concisely, we construct the sketch as
\begin{equation}\label{eq:YAQ}
\mathbf{Y} = \mathbf{A}\mathbf{\Omega},
\end{equation}
where $\mathbf{\Omega} \in \mathbb{R}^{m\times k}$ is a random test matrix with entries drawn from, for example, a Gaussian distribution.  
If $\mathbf{A}$ has exact rank $k$, then the sketch $\mathbf{Y}$ spans with high probability a basis for the column space of the source matrix. However, most data matrices do not feature exact rank (i.e., the singular values $\{\sigma_i\}_{i=k+1}^n$ are non-zero). Thus, instead of just using $k$ samples, it is favorable to construct a slightly oversampled sketch $\mathbf{Y} \in \mathbb{R}^{m\times l}$ which has $l=k+p$ columns, were $p$ denotes the number of additional samples. In most situations small values $p=\{10,20\}$ are sufficient to obtain a good basis that is comparable to the best possible basis~\citep{RandomizedMatrixComputations}. An orthonormal basis $\mathbf{Q} \in \mathbb{R}^{m\times l}$ is then obtained via the QR-decomposition  $\mathbf{Y}=\mathbf{Q}\mathbf{R}$, such that
\begin{displaymath}
\mathbf{A} \approx \mathbf{Q}\mathbf{Q}^\top\mathbf{A}
\end{displaymath}
is satisfied. Finally, the source matrix $\mathbf{A}$ is projected to low-dimensional space
\begin{equation}~\label{eq:BQA}
\mathbf{B} = \mathbf{Q}^\top\mathbf{A},
\end{equation}
where $\mathbf{B} \in \mathbb{R}^{l\times n}$. (Note, that Eq.~\ref{eq:BQA} requires a second pass over the source matrix.)
The matrix $\mathbf{B}$ can then be used to efficiently compute the matrix decomposition of interest, e.g., the SVD. The approximation error can be controlled by a combination of oversampling and power iteration~\citep{rokhlin2009randomized,halko2011rand,gu2015subspace}.  

Randomized matrix algorithms are not only pass efficient, but they also have the ability to exploit modern computational 
parallelized and distributed computing architectures.  Implementations in \textit{MATLAB}, \textit{C} and \textit{R} are provided by~\cite{szlam2014implementation},~\cite{voronin2015rsvdpack}, and~\cite{erichson2016randomized}.

\section{Randomized CP Decomposition} \label{sec:cCP}
Given a third order tensor $\tX \in \mathbb{R}^{I \times J \times K}$, the objective of the CP decomposition is to find a set of $R$ normalized rank-one tensors $\{\mathbf{a}_r \circ \mathbf{b}_r \circ \mathbf{c}_r\}_{r=1}^R$ which best approximates $\tX$, i.e., minimizes the Frobenius norm
\begin{equation}
\begin{aligned}
& \underset{ \hat{\tX} }{\text{minimize}}
& & \| \tX - \hat{\tX} \|_F^2 
& \text{subject to}
& & \hat{\tX} =\sum_{r=1}^{R} \mathbf{a}_r \circ \mathbf{b}_r \circ \mathbf{c}_r.
\end{aligned}
\label{eq:CPopt}
\end{equation}
The challenge is that this problem is highly nonconvex and unlike PCA there is no closed-form solution. Solution methods for this optimization problem therefore rely on iterative schemes (e.g., alternating least squares). These solvers are not guaranteed to find a global minimum, yet for many practical problems, they can find high-quality solutions. Iterative schemes, however, can be computationally demanding if the dimensions of $\tX$ are large.

Fortunately, only the column spaces of the modes are of importance for obtaining the factor matrices $\mathbf{A}$, $\mathbf{B}$, $\mathbf{C}$, rather then the individual columns of the mode matricizations $\tX_{(1)}$, $\tX_{(2)}$, $\tX_{(3)}$ of the tensor $\tX$. This is because the CP decomposition learns the components based on proportional variations in inter-point distances between the components. Therefore, a compressed tensor $\tB \in \mathbb{R}^{k \times k \times k}$ must preserve pairwise Euclidean distances, where $k\geq R$.
%
This in turn requires that column spaces, and thus pairwise distances, are approximately preserved by compression - this can be achieved by generalizing the concepts of randomized matrix algorithms to tensors. We build upon the methods introduced by \cite{Martinsson201147} and \cite{halko2011rand}, as well as related work on randomized tensors by~\cite{DrineasTensor}, who proposed a randomized algorithm based on random column selection. 

Figure~\ref{fig:rCParchitecture} shows the schematic of the randomized CP decomposition architecture. Our approach proceeds in two stages. The first stage, detailed in Section~\ref{sec:rand_tensor_alg}, applies random projections with power iteration to obtain a compressed tensor $\mathcal{B}$, with expressivity analysis in Section~\ref{sec:error_analysis}. In Section~\ref{sec:optimization}, we describe two algorithms for performing CP on the compressed tensor $\mathcal{B}$ and approximating the original CP factor matrices. Section~\ref{sec:implementation} provides additional details about our implementation in Python. 
\begin{figure}[!b]
	\begin{center}
		\begin{tikzpicture}[thick,scale=0.7, every node/.style={transform shape}]
		\node [blocky4,name=fdata] { $\tX$};       
		\node [blocky4,name=CP, right of=fdata, node distance=8cm]  {$[\![\bm{\bA,\bB,\bC}  ]\!]$};
		\node [blocky4,name=cdata,below of=fdata,node distance=5.0cm] {$\tB$};
		\node [blocky4,name=rCP,below of=CP,node distance=5.0cm] {$[\![\bm{ \tilde{\bA},\tilde{\bB},\tilde{\bC}}  ]\!]$};
		
		\path [draw, ->, line width=1.2pt] (fdata) -- node[name=toparrow,above of=fdata,node distance=.25cm]{\textbf{CP decomposition}}(CP);
		\path [draw, ->, line width=1.2pt, red] (cdata) -- node[name=bottomarrow,above of=cdata,node distance=.25cm]{\textbf{CP decomposition}}(rCP);
		\path [draw, ->, line width=1.2pt, red] (fdata) -- node[name=leftarrow, left of=fdata, node distance = 2.3cm, text width=4cm]{A degree of randomness is employed to derive a smaller tensor from the (big) input tensor.}(cdata);
		\path [draw, <-, line width=1.2pt, red] (CP) -- node[name=rightarrow, right of=CP, node distance = 2.0cm, text width=3cm]{Recover near-optimal factors.}(rCP);
		\node [name=data,above of=fdata,node distance=1.4cm] {\textbf{Tensor}};
		\node [name=modes,above of=CP,node distance=1.4cm] {\textbf{Factor matrices}};
		\node [name=full,left of=fdata,node distance=1.8cm]{\begin{sideways}\textbf{Big}\end{sideways}};
		\node [name=full,left of=cdata,node distance=1.8cm]{\begin{sideways}\textbf{Small}\end{sideways}};
		\end{tikzpicture}
	\end{center}
	\caption{Schematic of the randomized CP decomposition architecture. First, a degree of randomness is employed to derive a smaller tensor $\tB$ from the big tensor $\tX$. Then, the CP decomposition is performed on $\tB$. Finally, the near-optimal high-dimensional factor matrices $\bA, \bB$ and $\bC$ are recovered from the approximate factor matrices $\tilde{\bA}, \tilde{\bB}$ and $\tilde{\bC}$.}
	\label{fig:rCParchitecture}
\end{figure}

\subsection{Randomized Tensor Algorithm}\label{sec:rand_tensor_alg}
The aim is to use randomization as a computational resource to efficiently build a suitable basis that captures the action of the tensor $\tX$. Assuming an $N$-way tensor $\tX \in \mathbb{R}^{I_1 \times \cdots \times I_N}$, the aim is to obtain a smaller compressed tensor $\tB \in \mathbb{R}^{k \times \cdots \times k}$, so that its $N$ tensor modes capture the action of the input tensor modes. Hence, we seek a natural basis in the form of a set of orthonormal matrices $\{\mathbf{Q}_n \in \mathbb{R}^{I_n \times R_n}\}_{n=1}^{N}$, so that  
\begin{equation}
\tX \approx \tX \times_1 \bQ_1\bQ_1^\top \times_2 \cdots \times_N \bQ_N\bQ_N^\top .
\end{equation}
Here the operator $\times_n$ denotes tensor-matrix multiplication defined as follows.
\begin{Definition}
	The $n$-mode matrix product $\tX \times_n  \bQ_n\bQ_n^\top$ multiplies a tensor by the matrices $ \bQ_n\bQ_n^\top$ in mode $n$, i.e., each mode-$n$ fiber is multiplied by $ \bQ_n\bQ_n^\top$
	\begin{eqnarray*}
		\tM = \tX \times_n  \bQ_n\bQ_n^\top &   \Leftrightarrow  & \mathbf{M}_{(n)} = \bQ_n\bQ_n^\top\tX_{(n)}. 
	\end{eqnarray*}	
\end{Definition}
Given a fixed target rank $k$, these basis matrices can be efficiently obtained using a randomized algorithm. {\color{magenta} First}, the approximate basis for the $n$-th tensor mode is obtained by drawing $k$ random vectors 
$\mathbf{\omega}_1, \dots,\mathbf{\omega}_{k} \in \mathbb{R}^{\prod_{i\neq n}{I_i}}$ from a Gaussian distribution. (Note that if $N$ is large, it might be favorable to draw the entries of $\omega$ from a quasi random sequence, e.g., Halton or Sobol sequence. These so-called quasi-random numbers are known to be less correlated in high-dimensions.) Stacked together, the $k$ random vectors $\omega$ form the random test matrix $\mathbf{\Omega}_n \in \mathbb{R}^{\prod_{i\neq n}{I_i} \times k}$ used to sample the column space of $\tX_{(n)}\in \mathbb{R}^{I_n \times \prod_{i\neq n}{I_i}}$ as follows
\begin{equation}\label{eq:sampleMatrix}
\mathbf{Y}_n = \tX_{(n)}\mathbf{\Omega}_n,
\end{equation}
where $\mathbf{Y}_n \in \mathbb{R}^{I_n \times k}$ is the sketch.
The sketch serves as an approximate basis for the range of the $n$-th tensor mode. Probability theory guarantees that the set of random vectors $\{\mathbf{\omega}_i\}_{i=1}^{k}$ are linearly independent with high probability. Hence, the corresponding random projections $\mathbf{y}_1, \dots,\mathbf{y}_{k}$ efficiently sample the range of a rank deficient tensor mode $\tX_{(n)}$. The economic QR decomposition of the sketch $\mathbf{Y}_n = \bQ_n \mathbf{R}_n$ is then used to obtain a natural basis, so that $\mathbf{Q}_n \in \mathbb{R}^{I_n\times k}$ is orthonormal and has the same column space as $\mathbf{Y}_n$. The final step restricts the tensor mode to this low-dimensional subspace
\begin{eqnarray}
\mathbf{\tB}_n =  \tX \times_n  \mathbf{Q}_n^\top  &   \Leftrightarrow  &  \mathbf{B}_n =  \mathbf{Q}_n^\top \tX_{(n)}.
\end{eqnarray}	
Thus, after $N$ iterations a compressed tensor ${\tB}$ and a set of orthonormal matrices is obtained. Since this is an iterative algorithm, we set $\tX \gets \mathbf{\tB}_n$ after each iteration.  
The number of columns of the basis matrices form a trade-off between accuracy and computational performance. We aim to use as few columns as possible, yet allow an accurate approximation of the input tensor. Assuming that the tensor $\tX$ exhibits low-rank structure, or equivalently, the rank $R$ is much smaller than the ambient dimensions of the tensor, the basis matrices will be an efficient representation. {\color{black} In practice, the compression performance is improved through using oversampling, i.e., drawing $l=k+p$ random vectors where $k$ is the target rank and $p$ the oversampling parameter.} 

The randomized algorithm as presented requires that the mode-$n$ unfolding of the tensor has a rapidly decaying spectrum in order to achieve good performance. However, this assumption is often not suitable, and the spectrum exhibits slow decay if the tensor is compressed several times. To overcome this issue, the algorithm's performance can be substantially improved using power iterations~\citep{rokhlin2009randomized,halko2011rand,gu2015subspace}. Power iterations turn a slowly decaying spectrum into a rapidly decaying one by taking powers of the tensor modes. Thus, instead of sampling $\tX_{(n)}$ we sample from the following tensor mode
\begin{equation*}
\tX_{(n)}^q \coloneqq (\tX_{(n)}\tX_{(n)}^\top  )^q \tX_{(n)},
\end{equation*}
where $q$ denotes a small integer. This power operation enforces that the singular values $\sigma_j$ of $\tX_{(n)}^q$ are $\{\!\sigma_j^{2q+1}\!\}_j$. Instead of using \eqref{eq:sampleMatrix}, an improved sketch is computed
\begin{equation}\label{eq:IMPsampleMatrix}
\mathbf{Y}_n = (\tX_{(n)}\tX_{(n)}^\top  )^q \tX_{(n)} \mathbf{\Omega}_n.
\end{equation}
However, if \eqref{eq:IMPsampleMatrix} is implemented in this form the basis may be distorted due to round-off errors. Therefore in practice, normalized subspace iterations are used to form the sketch, meaning that the sketch is orthornormalized between each power iteration in order to stabilize the algorithm. For implementation, see~\cite{voronin2015rsvdpack} and~\cite{szlam2014implementation}.  

The combination of oversampling and additional power iterations can be used to control the trade-off between approximation quality and computational efficiency of the randomized tensor algorithm. Our results, for example, show that just $q=2$ subspace iterations and an oversampling parameter of about $p=10$ achieves near-optimal results.  Algorithm~\ref{alg:qb} summarizes the computational steps.
\begin{algorithm}[!t]
	\scalebox{0.95}{		
		\begin{minipage}{210mm}
			\begin{tabbing}
				\hspace{3mm} \= \hspace{5mm} \= \hspace{5mm} \= \hspace{5mm} \= \hspace{60mm} \=\kill
				\textbf{Require:} An $N$-way tensor $\mathbf{\tX}$, and a desired target rank $k$.\\[1mm]
				\textbf{Optional:} Parameters $p$ and $q$ to control oversampling, and power iterations.\\[1mm] 
				
				(1)  \> \> $\tB = \tX$ \> \> \>{\color{blue}$\textrm{initialize compressed tensor}$}\\[1mm]
				
				(2)  \> \> \textbf{for} $n = 1,\dots,N$ \> \> \>{\color{blue}$\textrm{iterate over all tensor modes}$} \\[1mm]				
				
				(3)  \> \> \>$l = k + p$ \> \> {\color{blue}$\textrm{slight oversampling}$}\\[1mm]
				
				(4)  \> \> \>$I,J = \texttt{dim}(\tB_{(n)})$ \> \> {\color{blue}$\textrm{dimension of the } n\textrm{-th tensor mode}$}\\[1mm]				
				
				(5)  \> \> \>$\mathbf{\Omega} = \texttt{rand}(J,l)$ \> \> {\color{blue}$\textrm{generate random test matrix}$}\\[1mm]
				
				(6)  \> \> \>$\mathbf{Y} = \tB_{(n)} \mathbf{\Omega}$ \> \> {\color{blue}$\textrm{compute sampling matrix}$}\\[1mm]
				
				(7)  \> \> \>\textbf{for} $j = 1,\dots,q$ \> \> {\color{blue}$\textrm{normalized power iterations (optional)}$} \\[1mm]
				
				(8)  \> \> \> \>$[\mathbf{Q}, \sim] = \texttt{lu}(\mathbf{Y})$ \> \\[1mm] 
				
				(9)  \> \> \> \>$[\mathbf{Z}, \sim] = \texttt{lu}(\tB_{(n)}^\top \mathbf{Q})$ \\[1mm]
				
				(10)  \> \> \> \>$\mathbf{Y} = \tB_{(n)} \mathbf{Z}$ \\[1mm]
				
				(11)  \> \> \>\textbf{end for}\\[1mm]
				
				(12)  \> \> \>$[\mathbf{Q}_n, \sim] = \texttt{qr}(\mathbf{Y})$ \> \> {\color{blue}\textrm{orthonormalize sampling matrix}}
				\\[1mm]
				
				(13)  \> \> \>$\mathbf{\tB} =  \tB \times_n  \mathbf{Q}_n^\top$ \> \> {\color{blue}$\textrm{project tensor to smaller space}$} \\[1mm]			
				
				(14)  \> \> \textbf{end for}\\[1mm]								
				
				\textbf{Return:} Compressed tensor $\tB$ of dimension $l\times \cdots \times l$, and a set of orthonormal \\ basis matrices $\{\mathbf{Q}_n \in \mathbb{R}^{I_n \times l}\}_{n=1}^N$.			
			\end{tabbing}
	\end{minipage}}	
	\begin{remark}  Due to the convenient mathematical properties of the normal distribution, a Gaussian random test matrix is {\color{black} often assumed in theoretical results}, however, in practice a uniform distributed random test matrix is sufficient. The performance can be further improved using structured random matrices~\citep{woolfe2008fast}. If information is uniformly distributed across the data, randomly selected columns can also be used to build a suitable basis as well, which avoids the matrix multiplication in step (6).
	\end{remark}		
	\begin{remark} For numerical stability, normalized power iterations using the pivoted LU decomposition are computed in step 7-11. We recommend a default value of $q=2$. 
	\end{remark}
	\begin{remark} In practice, the user can decide which modes to compress and specify different oversampling parameters for these modes. 
	\end{remark}				
	\caption{A prototype randomized tensor compression algorithm.}
	\label{alg:qb}
\end{algorithm}

\subsubsection{Expressivity analysis}\label{sec:error_analysis}
The average behavior of the randomized tensor algorithm is characterized using the expected residual error  
\begin{equation}
\E\|\tE \|_F = \| \tX - \hat{\tX} \|_F,
\end{equation}
where $\hat{\tX} = \tX \times_1 \bQ_1\bQ_1^\top \times_2 \cdots \times_N \bQ_N\bQ_N^\top$. 
Theorem~\ref{thm:ubound} {\color{black} generalizes the matrix version of Theorem 10.5 formulated by~\cite{halko2011rand} to the tensor case}.  
\begin{Theorem}[Expected Frobenius error]\label{thm:ubound} 
	Consider a low-rank real $N$-way tensor $\tX \in \mathbb{R}^{I_1 \times \cdots \times I_N}$. Then the expected approximation error, given a target rank $k\geq2$ and an oversampling parameter $p\geq2$ for each mode, is 
	\begin{eqnarray}
	\E\| \tX - \tX \times_1 \bQ_1\bQ_1^\top \times_2 \cdots \times_N \bQ_N\bQ_N^\top \|_F & \leq & \sqrt{ 1 + \frac{k}{p-1}} \cdot \sqrt{\sum_{n=1}^{N}\sum_{j>k}\sigma_{nj}^2},
	\label{eq:theorem}
	\end{eqnarray}
	where $\sigma_{nj}$ denotes the $j$ singular value of the mode-$n$ unfolding of the source tensor $\tX$.
\end{Theorem}
For the proof see Appendix~\ref{app:proof1}. Intuitively, the projection of each tensor mode onto a low-dimensional space introduces an additional residual. This is expressed by the double sum on the right hand side. If the low-rank approximation captures the column space of each mode accurately, then the singular values $j>k$ for each mode $n$ are small. Moreover, the error can be improved using the oversampling parameter. The computation of additional power (subspace) iterations can improve the error further. This result again follows by generalizing the results of~\cite{halko2011rand} to tensors. Sharper performance bounds for both oversampling and additional power iterations can be derived following, for instance, the results by~\cite{witten2015randomized}.

\subsection{Optimization Strategies}\label{sec:optimization}
{\color{black} There are several optimization strategies for minimizing the objective function defined in \eqref{eq:cCPopt}, of which} we consider, alternating least squares (ALS) and block coordinate descent (BCD). Both methods are suitable to operate on a compressed tensor $\tB \in \mathbb{R}^{k \times \cdots \times k}$, where $k\geq R$. The optimization problem \eqref{eq:CPopt} is reformulated 
\begin{equation}\label{eq:cCPopt}
\begin{aligned}
& \underset{ \hat{\tB} }{\text{minimize}}
& & \| \tB - \hat{\tB} \|_F^2 
& \text{subject to}
& & \hat{\tB} =\sum_{r=1}^{R} \mathbf{\tilde{a}}_r \circ \mathbf{\tilde{b}}_r \circ \mathbf{\tilde{c}}_r.
\end{aligned}
\end{equation}
Once the compressed factor matrices $\mathbf{\tilde{A}} \in \mathbb{R}^{k \times R}$, $\mathbf{\tilde{B}} \in \mathbb{R}^{k \times R}$, $\mathbf{\tilde{C}} \in \mathbb{R}^{k \times R}$ are estimated, the full factor matrices can be recovered 
\begin{eqnarray}\label{eq:qbs}
\mathbf{A}  \approx \mathbf{Q}_1\mathbf{\tilde{A}}, \,\, 
\mathbf{B}  \approx  \mathbf{Q}_2\mathbf{\tilde{B}}, \,\,
\mathbf{C}  \approx  \mathbf{Q}_3\mathbf{\tilde{C}}, 
\end{eqnarray}
where $\mathbf{Q}_1 \! \in \! \mathbb{R}^{I \times k}$, $\mathbf{Q}_2 \! \in \! \mathbb{R}^{J \times k}$, $\mathbf{Q}_3 \! \in \! \mathbb{R}^{K \times k}$ denote the orthonormal basis matrices. For simplicity we focus on third order tensors, but the result generalizes to $N$-way tensors.

\subsubsection{ALS Algorithm}
Due to its simplicity and efficiency, ALS is the most popular method for computing the CP decomposition~\citep{comon2009tensor,kolda2009tensor}. We note that the optimization \eqref{eq:cCPopt} is equivalent to
\begin{equation*}
\begin{aligned}
\underset{\bf {\tilde{A}, \tilde{B}, \tilde{C}} }{\text{minimize}} & & \| \tB - \sum_{r=1}^{R} \mathbf{\tilde{a}}_r \circ \mathbf{\tilde{b}}_r \circ \mathbf{\tilde{c}}_r \|_F^2 
\end{aligned}
\end{equation*}
with respect to the factor matrices $\bf {\tilde{A}, \tilde{B}}$ and $\bf {\tilde{C}}$. The tensor $\tB$ can further be expressed in matricized form
\begin{eqnarray*}
	\mathbf{\tB}_{(1)}  \approx  \mathbf{\tilde{A}}(\mathbf{\tilde{C}}\odot\mathbf{\tilde{B}})^{\top}\!, \,\,
	\mathbf{\tB}_{(2)}  \approx  \mathbf{\tilde{B}}(\mathbf{\tilde{C}}\odot\mathbf{\tilde{A}})^{\top}\!, \,\,
	\mathbf{\tB}_{(3)}  \approx  \mathbf{\tilde{C}}(\mathbf{\tilde{B}}\odot\mathbf{\tilde{A}})^{\top}\!,
\end{eqnarray*}
where $\odot$ denotes the Khatri-Rao product. The optimization problem in this form is non-convex. However, an estimate for the factor matrices can be obtained using the least-squares method {\color{black} as follows.}  
The ALS algorithm updates one component, while holding the other two components fixed, in an alternating fashion until convergence. It iterates over the following subproblems
\begin{eqnarray}
\mathbf{\tilde{A}}^{j+1} & = & \argmin_{\mathbf{\tilde{A}}}\|\tB_{(1)} - \mathbf{\tilde{A}}(\mathbf{\tilde{C}}^{j}\odot\mathbf{\tilde{B}}^{j})^{\top}  \|,
\label{eq:Astar} \\
\mathbf{\tilde{B}}^{j+1} & = & \argmin_{\mathbf{\tilde{B}}}\|\tB_{(2)} - \mathbf{\tilde{B}}(\mathbf{\tilde{C}}^{j}\odot\mathbf{\tilde{A}}^{j+1})^{\top}  \|, \\
\mathbf{\tilde{C}}^{j+1} & = & \argmin_{\mathbf{\tilde{C}}}\|\tB_{(3)} - \mathbf{\tilde{C}}(\mathbf{\tilde{B}}^{j+1}\odot\mathbf{\tilde{A}}^{j+1})^{\top}  \|.
\end{eqnarray}
Each step therefore involves a least-squares problem which can be solved using the Khatri-Rao product pseudo-inverse. 
Algorithm~\ref{alg:als} summarizes the computational steps. 
\begin{Definition} The Khatri-Rao product pseudo-inverse is defined as
	\begin{equation*}
	(\bA \odot \bB)^{\dagger}=(\bA^\top \bA * \bB^\top \bB)^{\dagger}(\bA \odot \bB)^{\top},
	\end{equation*}
	where the operator $*$ denotes the Hadamard product, i.e., the elementwise multiplication of two equal sized matrices.
\end{Definition}
%

There exist few general convergence guarantees for the ALS algorithm~\citep{uschmajew2012local,wang2014global}. Moreover, the final solution tends to depend on the initial guess $\bf {\tilde{A}}^0,\bf{\tilde{B}}^0$ and $\bf {\tilde{C}}^0$. A standard initial guess uses the eigenvectors of $\tB_{(1)}\tB_{(1)}^\top$, $\tB_{(2)}\tB_{(2)}^\top$, $\tB_{(3)}\tB_{(3)}^\top$~\citep{TTB_Software}. Further, it is important to note that normalization of the factor matrices is necessary after each iteration to achieve good convergence. This prevents singularities of the Khatri-Rao product pseudo-inverse~\cite{kolda2009tensor}. The algorithm can be further improved by reformulating the above subproblems as regularized least-squares problems, see for instance  \cite{li2013some} for technical details and convergence results. {\color{black} The structure imposed by the ALS algorithm on the factor matrices} permits the formulation of non-negative, or sparsity-constrained tensor decompositions~\citep{cichocki2009nonnegative}.  
\begin{algorithm}[t]
	\scalebox{0.95}{		
		\begin{minipage}{210mm}
			\begin{tabbing}
				\hspace{3mm} \= \hspace{5mm} \= \hspace{5mm} \= \hspace{65mm} \=\kill
				\textbf{Require:} An $I\times J\times K$ tensor $\mathbf{\tX}$, and a desired target rank $R$.\\[1mm]
				\textbf{Optional:} Parameters $p$ and $q$ to control oversampling, and power iterations.\\[1mm]

				(1)  \> \> $\tB ,\mathbf{Q}_A, \mathbf{Q}_B , \mathbf{Q}_C  = \texttt{compress}(\tX,R,p,q)$ \> \> {\color{blue}$\textrm{compress tensor using Algorithm~\ref{alg:qb}}$}\\[1mm]
				
				(2)  \> \> $\mathbf{B,C} = [\texttt{eig}(\tB_{(2)},\tB_{(2)}^\top), \texttt{eig}(\tB_{(3)},\tB_{(3)}^\top)]$ \> \> {\color{blue}$\textrm{use first } R \textrm{ eigenvectors for initialization}$}\\[1mm]

				(3)  \> \> \textbf{repeat}  \> \>  \\[1mm]
				
				(4)  \> \> \> $\bA = \tB_{(1)} (\bC \odot \bB) (\bC^\top \bC * \bB^\top \bB)^\dagger$ \> \\[1mm]  
				
				(5)  \> \> \> $\bA = \bA / \texttt{norm}(\bA)$ \> \\[1mm]  				
				
				(6)  \> \> \> $\bB = \tB_{(2)} (\bC \odot \bA) (\bC^\top \bC * \bA^\top \bA)^\dagger$ \> \\[1mm] 
				
				(7)  \> \> \> $\bB = \bB / \texttt{norm}(\bB)$ \> \\[1mm]  								 
				
				(8)  \> \> \> $\bC = \tB_{(3)} (\bB \odot \bA) (\bB^\top \bB * \bA^\top \bA)^\dagger$ \> \\[1mm]  
				
				(9)  \> \> \> $\bm{\lambda} = \texttt{norm}(\bC)$ \> \\[1mm]  										
				
				(10)  \> \> \> $\bC = \bC / \bm{\lambda}$ \> \\[1mm]  		
				
				(11)  \> \> \textbf{until} {convergence criterion is reached}\\[1mm]

				(12)  \> \> $\bA,\bB,\bC = [\mathbf{Q}_A \bA, \mathbf{Q}_B \bB, \mathbf{Q}_C \bC]$ \> \> {\color{blue}$\textrm{recover factor matrices}$} \\[1mm]
				
				(13)  \> \> {re-normalize the factor matrices and update the scaling vector $\bm{\lambda}$}\\[1mm]		
				
				\textbf{Return:} Normalized factor matrices $\bA,\bB,\bC$ and the scaling vector $\bm{\lambda}$.			
			\end{tabbing}
	\end{minipage}}	
	\caption{A prototype randomized CP algorithm using ALS.}
	\label{alg:als}
\end{algorithm}

\subsubsection{BCD Algorithm}
While ALS is the most popular algorithm for computing the CP decomposition, many alternative algorithms have been developed. One {\color{black} alternate approach is based on block coordinate descent,} BCD~\citep{BlockCoordinateDescent}. \cite{cichocki2009fast} first proposed this approach for computing nonnegative tensor factorizations. 
The BCD algorithm is based on the idea of successive rank-one deflation. Unlike ALS, which updates the entire factor matrix at each step, BCD computes the rank-1 tensors in a hierarchical fashion. Therefore, the algorithm treats each component $\ba_r,\bb_r,\bc_r$  as a block. First, the most correlated rank-1 tensor is computed; then the second most correlated rank-1 tensor is learned on the residual tensor, and so on. Assuming that $\tilde{R}=r-1$ components have been computed, then the $r$-th compressed residual tensor $\tY_{\text{res}}$ is defined 
\begin{equation}
\tY_{\text{res}} = \tB - \sum_{r=1}^{\tilde{R}} \mathbf{\tilde{a}}_r \circ \mathbf{\tilde{b}}_r \circ \mathbf{\tilde{c}}_r.
\end{equation}
The algorithm then iterates over the following subproblems
\begin{eqnarray}
\mathbf{\tilde{a}}_r^{j+1} & = & \argmin_{\mathbf{\tilde{a}}_r}\|\tY_{\text{res}(1)} - \mathbf{\tilde{a}}_r(\mathbf{\tilde{c}}_r^{j}\odot\mathbf{\tilde{b}}_r^{j})^{\top}  \|, \\
\mathbf{\tilde{b}}_r^{j+1} & = & \argmin_{\mathbf{\tilde{b}_r}}\|\tY_{\text{res}(2)} - \mathbf{\tilde{b}}_r(\mathbf{\tilde{c}}_r^{j}\odot\mathbf{\tilde{a}}_r^{j+1})^{\top}  \|, \\
\mathbf{\tilde{c}}_r^{j+1} & = & \argmin_{\mathbf{\tilde{c}}_r}\|\tY_{\text{res}(3)} - \mathbf{\tilde{c}}_r(\mathbf{\tilde{b}}_r^{j+1}\odot\mathbf{\tilde{a}}_r^{j+1})^{\top}  \|.
\end{eqnarray}
Note that the computation can be more efficiently evaluated without explicitly constructing the residual tensor $\tY_\text{res}$~\citep{kim2014algorithms}. Algorithm~\ref{alg:bcd} summarizes the computation.
\begin{algorithm}[t]
	\scalebox{0.95}{		
		\begin{minipage}{210mm}
			\begin{tabbing}
				\hspace{3mm} \= \hspace{5mm} \= \hspace{5mm} \= \hspace{5mm} \= \hspace{60mm} \=\kill
				\textbf{Require:} An $I\times J\times K$ tensor $\mathbf{\tX}$, and a desired target rank $R$.\\[1mm]
				\textbf{Optional:} Parameters $p$ and $q$ to control oversampling, and power iterations.\\[1mm]

				(1)  \> \> $\tB ,\mathbf{Q}_A, \mathbf{Q}_B , \mathbf{Q}_C  = \texttt{compress}(\tX,R,p,q)$ \> \> \>{\color{blue}$\textrm{compress tensor using Algorithm~\ref{alg:qb}}$}\\[1mm]
				
				(2)  \> \> $\mathbf{B,C} = [\texttt{eig}(\tB_{(2)},\tB_{(2)}^\top), \texttt{eig}(\tB_{(3)},\tB_{(3)}^\top)]$ \> \> \> {\color{blue}$\textrm{use first } R \textrm{ eigenvectors for initialization}$}\\[1mm]
				
				(3)  \> \> $\tY = \tB$ \> \> \>{\color{blue}$\textrm{initialize residual tensor}$} \\[1mm]			
				
				(4)  \> \> \textbf{for} $r = 1,\dots,R$ \> \> \>{\color{blue}$\textrm{compute rank-$r$ approximation}$} \\[1mm]
				
				(5)  \> \> \> \textbf{repeat}  \> \>  \\[1mm]
				
				(6)  \> \> \> \> $\ba_r = \tY_{(1)} (\bc_r \odot \bb_r) (\bc_r^\top \bc_r * \bb_r^\top \bb_r)^\dagger$ \> \\[1mm]  
				
				(7)  \> \> \> \> $\ba_r = \ba_r / \texttt{norm}(\ba_r)$ \> \\[1mm]  				
				
				(8)  \> \> \> \>$\bb_r = \tY_{(2)} (\bc_r \odot \ba_r) (\bc_r^\top \bc_r * \ba_r^\top \ba_r)^\dagger$ \> \\[1mm] 
				
				(9)  \> \> \> \>$\bb_r = \bb_r / \texttt{norm}(\bb_r)$ \> \\[1mm]  								 
				
				(10)  \> \> \> \>$\bc_r = \tY_{(3)} (\bb_r \odot \ba_r) (\bb_r^\top \bb_r * \ba_r^\top \ba_r)^\dagger$ \> \\[1mm]  
				
				(11)  \> \> \> \>$\bm{\lambda}_r = \texttt{norm}(\bc_r)$ \> \\[1mm]  										
				
				(12)  \> \> \> \>$\bc_r = \bc_r / \bm{\lambda}_r$ \> \\[1mm]  
				
				(13)  \> \> \> \textbf{until} {convergence criterion is reached}\\[1mm]	
				
				(14)  \> \> \> $\tY = \tB - [\![\bm{\lambda}_{[1:r]};\bA_{[:,1:r]},\bB_{[:,1:r]},\bC_{[:,1:r]}]\!]$ \> \> {\color{blue}$\textrm{update residual tensor}$} \\[1mm]																
				
				(15)  \> \> \textbf{end for}\\[1mm]				
				
				(16)  \> \> $\bA,\bB,\bC = [\mathbf{Q}_A \bA, \mathbf{Q}_B \bB, \mathbf{Q}_C \bC]$ \> \> \>{\color{blue}$\textrm{recover factor matrices}$} \\[1mm]
				
				(17)  \> \> {re-normalize the factor matrices and update the scaling vector $\bm{\lambda}$}\\[1mm]		
				
				\textbf{Return:} Normalized factor matrices $\bA,\bB,\bC$ and the scaling vector $\bm{\lambda}$.			
			\end{tabbing}
	\end{minipage}}	
	\caption{A prototype randomized CP algorithm using BCD.}
	\label{alg:bcd}
\end{algorithm}

\subsection{Implementation Details}\label{sec:implementation}
The algorithms we present are implemented in the programming language \textit{Python}, using numerical linear algebra tools provided by the \textit{SciPy} (Open Source Library of Scientific Tools) package~\citep{scipy}. \textit{SciPy} provides MKL (Math Kernel Library) accelerated high performance implementations of \textit{BLAS} and \textit{LAPACK} routines. Thus, all linear algebra operations are threaded and highly optimized on Intel processors. 
%
The implementation of the CP decomposition follows the \textit{MATLAB Tensor Toolbox} implementation~\citep{TTB_Software}. This implementation normalizes the components after each step to achieve better convergence. Furthermore, we use eigenvectors (see above) to initialize the factor matrices. Interestingly, randomly initialized factor matrices have the ability to achieve slightly better approximation errors, but re-running the algorithms several times with different random seeds can display significant variance in the results. Hence only the former approach is used for initialization. We note that the randomized algorithm introduces some randomness and slight variations into the CP decompositions as well. However, randomization can also act as an implicit regularization on the CP decomposition~\citep{Mahoney2011}, meaning that the results of the randomized algorithm can be in some cases even `better' than the results of the corresponding deterministic implementation. {\color{black} Finally, the convergence criterion is defined as the change in fit, following~\cite{TTB_Software}}. The algorithm therefore stops when the improvement of the fit $\rho$ is less then a predefined threshold, where the fit is computed using 
%
$\rho = 1 - ({\|\tX\|_F^2 + \| \hat{\tX} \|_F^2 - 2\cdot \langle \hat{\tX},\tX\rangle})/{\|\tX\|_F^2 }$.

\section{Numerical Results} \label{sec:results}

The randomized CP algorithm is evaluated on a number of examples where the near optimal approximation of massive tensors can be achieved in a fraction of the time using the randomized algorithm. 
Approximation accuracy is computed with the relative error 
${\|\tX - \hat{\tX}\|_F}/{\|\tX\|_F }$,
where $\hat{\tX}$ denotes the approximated tensor. 



%
\begin{figure}[!b]
	\centering
	\DeclareGraphicsExtensions{.pdf}
	\includegraphics[width=0.85\textwidth]{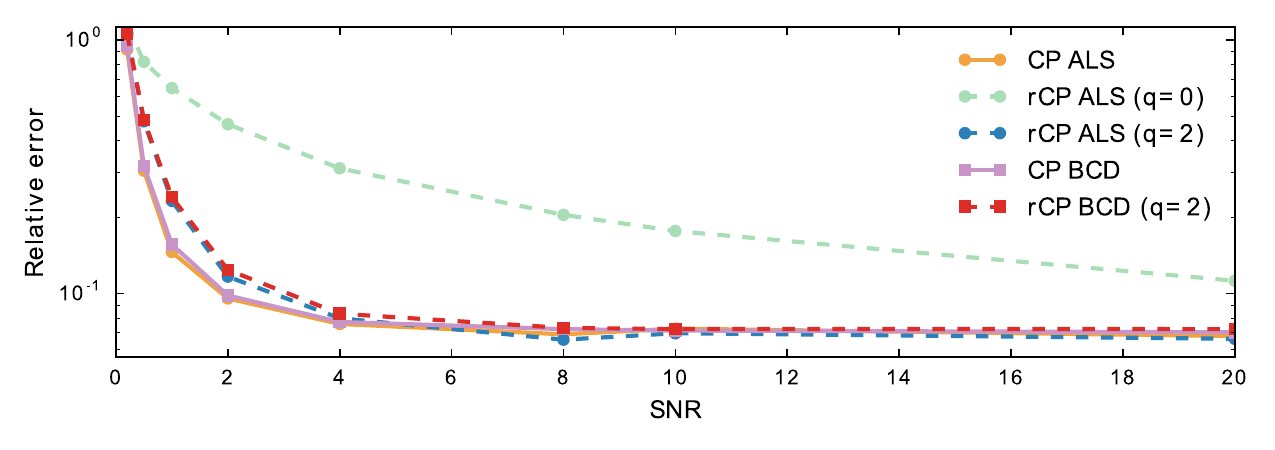}
	\vspace{-0.1in}		
	\caption{Average relative error, plotted on a log scale, against increasing signal to noise ratio. The analysis is performed on a rank $R=50$ tensor of dimension $100\times 100\times 100$. Power iterations improve the the approximation accuracy considerably, while the performance without power iterations is poor.}
	\label{fig:snr}
\end{figure}

\subsection{Computational Performance}\label{sec:compp}

\begin{figure}[!b]
	\centering
	\begin{subfigure}[t]{0.85\textwidth}
		\centering
		\DeclareGraphicsExtensions{.pdf}
		\includegraphics[width=1\textwidth]{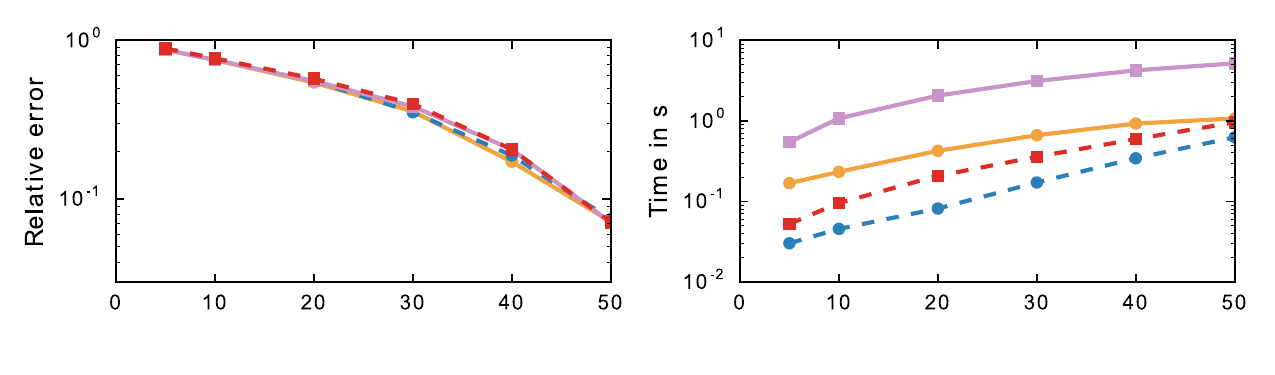}
		\vspace{-0.45in}
		\caption{Tensor of dimension $100\times 100\times 100$.}
	\end{subfigure}
	\begin{subfigure}[t]{0.85\textwidth}
		\centering
		\DeclareGraphicsExtensions{.pdf}
		\includegraphics[width=1\textwidth]{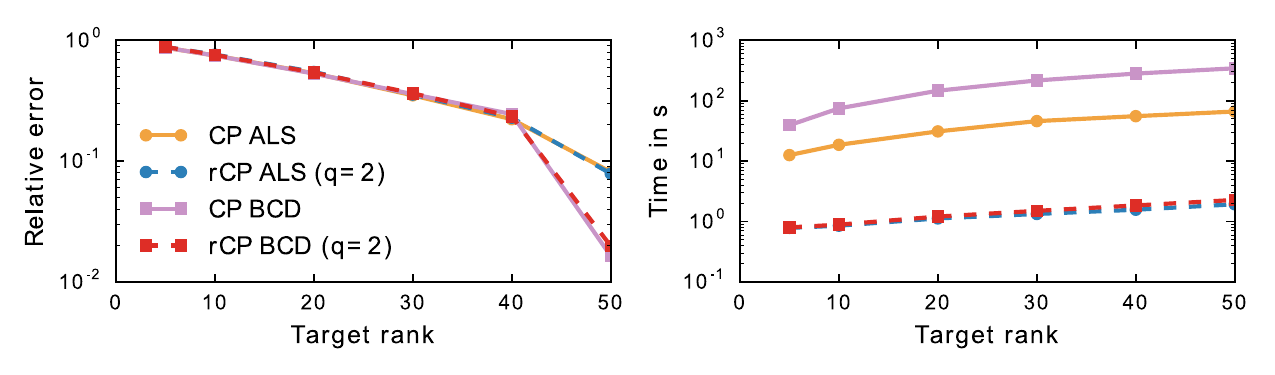}
		\vspace{-0.25in}
		\caption{Tensor of dimension $400\times 400\times 400$.}
	\end{subfigure}	
	
	\begin{subfigure}[t]{0.85\textwidth}
		\centering
		\DeclareGraphicsExtensions{.pdf}
		\includegraphics[width=1\textwidth]{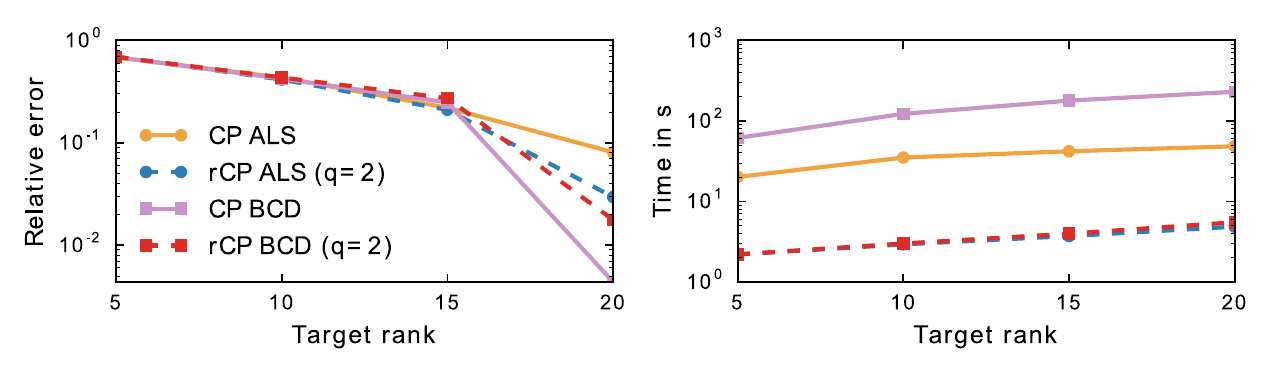}
		\caption{Tensor of dimension $100\times 100\times 100\times 100$.  }
		\label{fig:timing_4dim}
	\end{subfigure}
	
	\caption{Random tensor approximation and performance for rank $R=50$ tensors:  rCP methods achieve speedups by $1$-$2$ orders of magnitude and the same accuracy as their deterministic counterpart.}
	\label{fig:performance}
\end{figure}
%

The robustness of the randomized CP algorithm is first assessed on random low-rank tensors. Here we illustrate the approximation quality in the presence of additive white noise. Figure~\ref{fig:snr} shows the average relative errors over $100$ runs for varying signal-to-noise ratios (SNR). 
In the case of high SNRs, all algorithms converge towards the same relative error. However, at excessive levels of noise (i.e., SNR$<4$) the deterministic CP algorithms exhibit small gains in accuracy over the randomized algorithms using $q=2$ power iterations, with which both the ALS and BCD algorithm show the same performance. The performance of the randomized algorithm without power iterations ($q=0$) is, however, poor, and stresses the importance of the power operation for real applications. The oversampling parameter for the randomized algorithms is set to $p=10$. Increasing $p$ can further improve accuracy. 

Next, the reconstruction errors and runtimes for tensors of varying dimensions are compared. Figure~\ref{fig:performance} shows the average evaluation results over $100$ runs for random low-rank tensors of different dimensions, and for varying target ranks $k$. The randomized algorithms achieve near optimal approximation accuracy while demonstrating substantial computational savings.
Interestingly, the relative error achieved by the BCD decreases sharply by about one order of magnitude when the target rank $k$ matches the actual tensor rank (here $R=50$). 

The computational advantage becomes more pronounced with increasing tensor dimensions, and as the number of iterations required for convergence increase. Using random tensors as presented here, all algorithms rapidly converge after about $4$ to $6$ iterations. However, it is evident that the computational cost per iteration of the randomized algorithm is substantially lower. Thus, the computational savings can be even more substantial in real world applications that may require several hundred iterations to converge.
Overall, the ALS algorithm is computationally more efficient than BCD, i.e., the deterministic ALS algorithm is faster than the BCD by nearly one order of magnitude, while the randomized algorithms exhibit similar computational timings.

Similar performance results are achieved for higher order tensors. 
Figure~\ref{fig:timing_4dim} shows the computational performance for a $4$-way tensor of dimension $100\times 100\times 100\times 100$. Again, the randomized algorithms achieve speedups of 1-2 orders of magnitude, while attaining good approximation errors. 

\subsection{Examples}\label{sec:realexamples}

The examples demonstrate the advantages of the randomized CP decomposition. The first is a multiscale toy video example, and the second is a simulated flow field behind a stationary cylinder. Due to the better and more natural interpretability of the BCD algorithm, only this algorithm is considered in subsequent sections.

\begin{figure}[!b]
	\centering
	\DeclareGraphicsExtensions{.pdf}
	
	\begin{overpic}[width=0.98\textwidth]{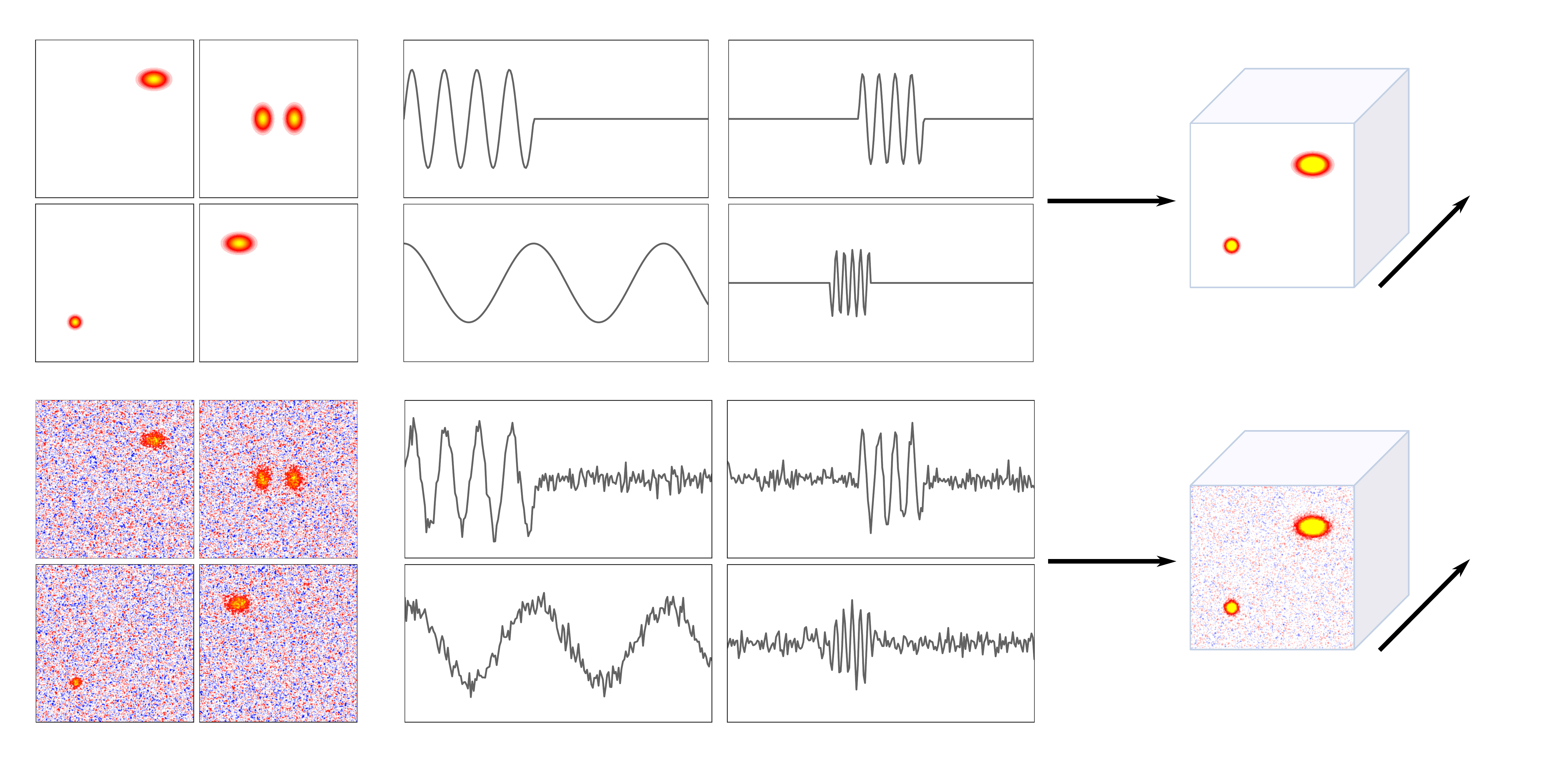} 
		\put(9,48){\small modes}
		\put(38,48){\small time dynamics}
		\put(80,48){\small tensor}
		\put(92,30){\small time}	
		\put(92,7){\small time}			
	\end{overpic}\vspace{-0.1cm}		
	
	\caption{Illustration of the multiscale toy video. The system is governed by four spatial modes experiencing intermittent oscillations in the temporal direction. The top row shows the clean system, while the the bottom row shows the noisy system with a signal-to-noise ratio of $2$. }
	\label{fig:toy_video}
\end{figure}

\subsubsection{Multiscale Toy Video Example}
The approximation of the underlying spatial modes and temporal dynamics of a system is a common problem in signal processing. In the following, we consider a toy example presenting multiscale intermittent dynamics in the time direction. The data consists of four Gaussian modes, each undergoing different frequencies of intermittent oscillation, for $215$ times steps in the temporal direction on a two-dimensional spatial grid  ($200\times 200$). The resulting tensor is of dimension $200\times 200\times 215$. Figure~\ref{fig:toy_video} shows the corresponding modes and the time dynamics. 
This problem becomes even more challenging when the underlying structure needs to be reconstructed from noisy measurements. 
Traditional matrix decomposition techniques such as the SVD are known to  face difficulties capturing such intermittent, multiscale structure.

\begin{figure}[!t]
	\centering
	
	\begin{subfigure}[t]{0.30\textwidth}
		\centering
		\DeclareGraphicsExtensions{.pdf}
		
		\begin{overpic}[width=1\textwidth]{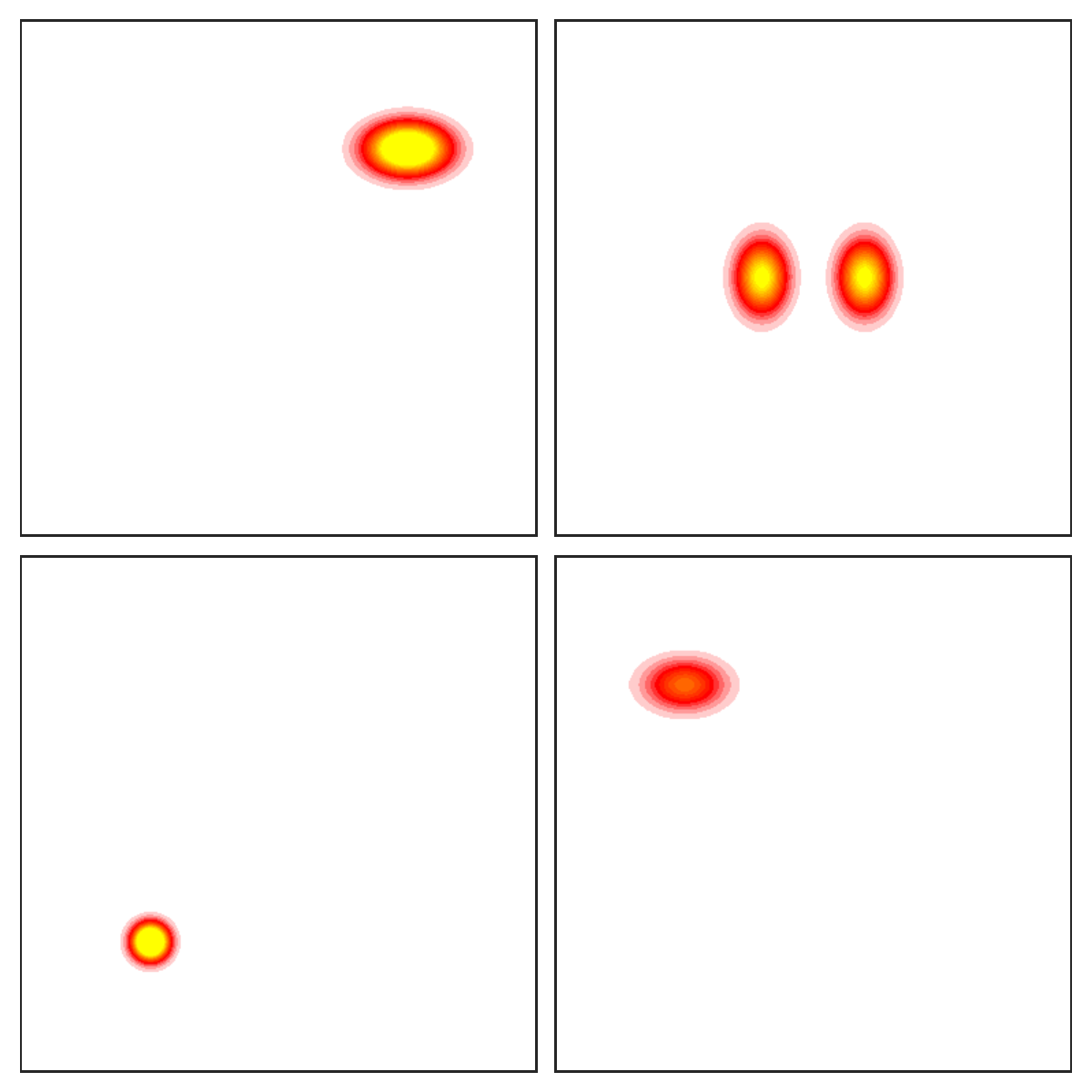}
			\put(-10,38){\rotatebox{90}{\small modes}}
		\end{overpic}
		
	\end{subfigure}	
	~
	\begin{subfigure}[t]{0.30\textwidth}
		\centering
		\DeclareGraphicsExtensions{.pdf}
		\includegraphics[width=1\textwidth]{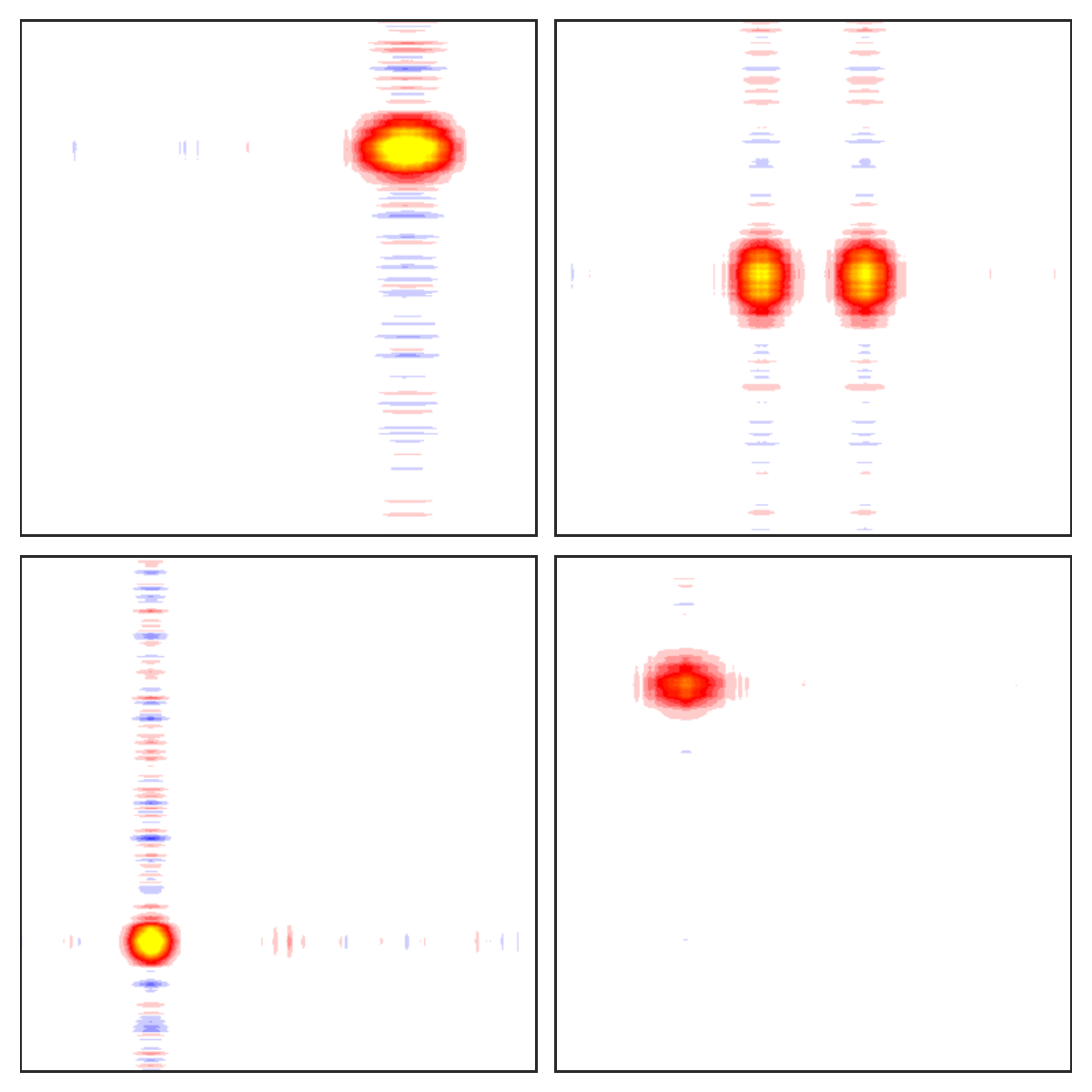}
		
	\end{subfigure}	
	~
	\begin{subfigure}[t]{0.30\textwidth}
		\centering
		\DeclareGraphicsExtensions{.pdf}
		\includegraphics[width=1\textwidth]{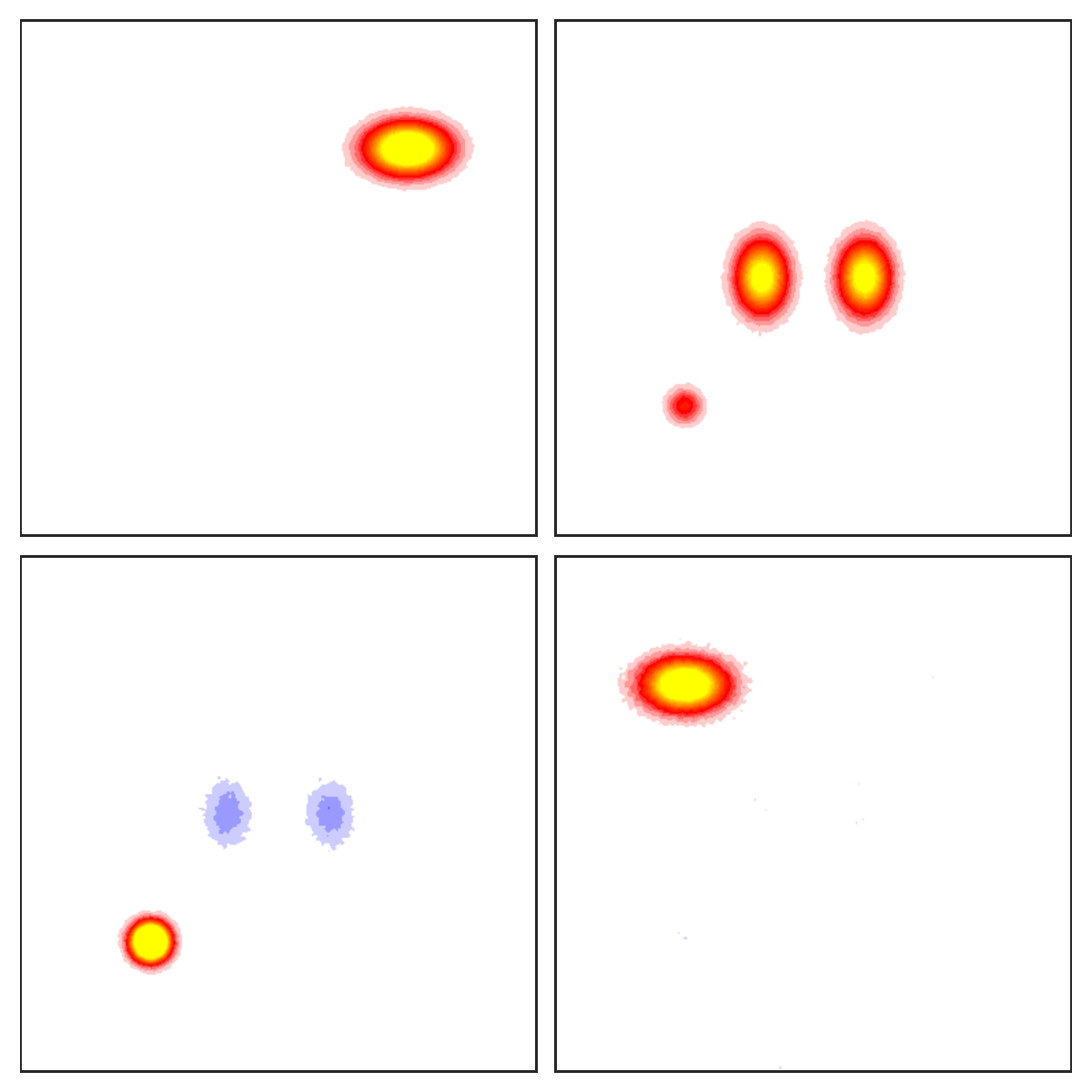}
		
	\end{subfigure}

	\begin{subfigure}[t]{0.30\textwidth}
		\centering
		\DeclareGraphicsExtensions{.pdf}
		\begin{overpic}[width=1\textwidth]{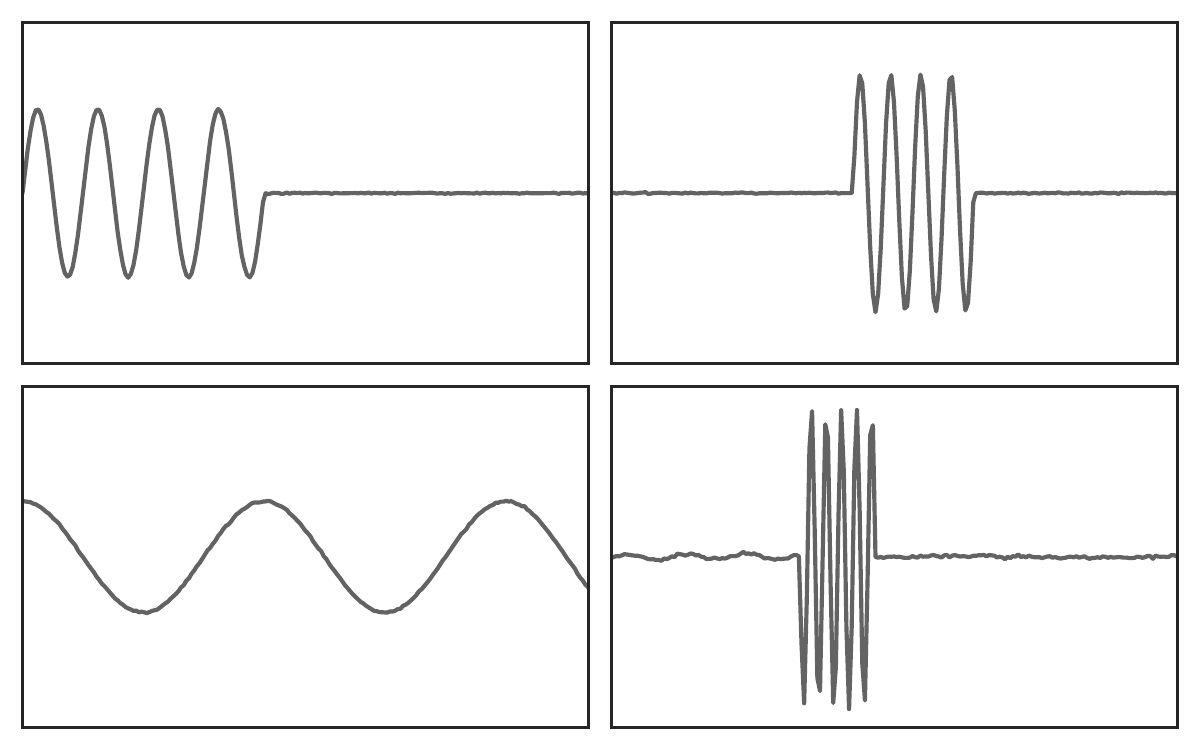}
			\put(-10,5){\rotatebox{90}{\small time dynamics}}
		\end{overpic}
		\caption{rCP ($q=2$).}
	\end{subfigure}	
	~
	\begin{subfigure}[t]{0.30\textwidth}
		\centering
		\DeclareGraphicsExtensions{.pdf}
		\includegraphics[width=1\textwidth]{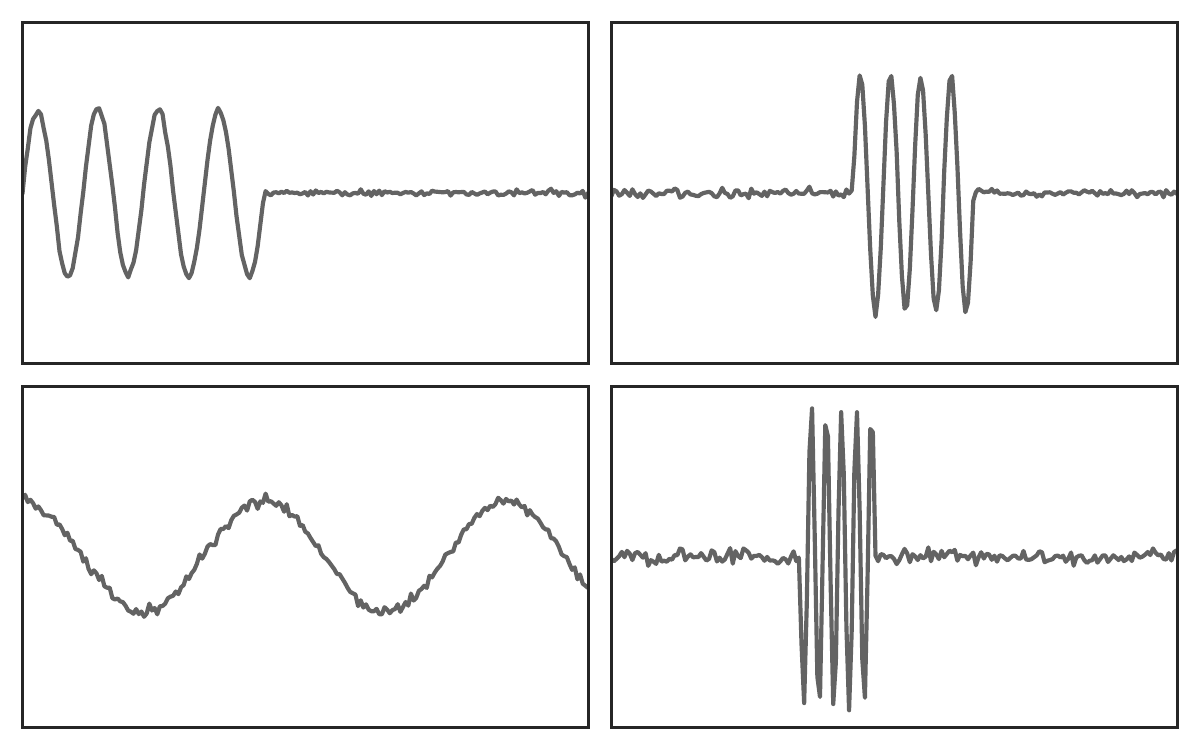}
		\caption{rCP ($q=0$).}
	\end{subfigure}	
	~
	\begin{subfigure}[t]{0.30\textwidth}
		\centering
		\DeclareGraphicsExtensions{.pdf}
		\includegraphics[width=1\textwidth]{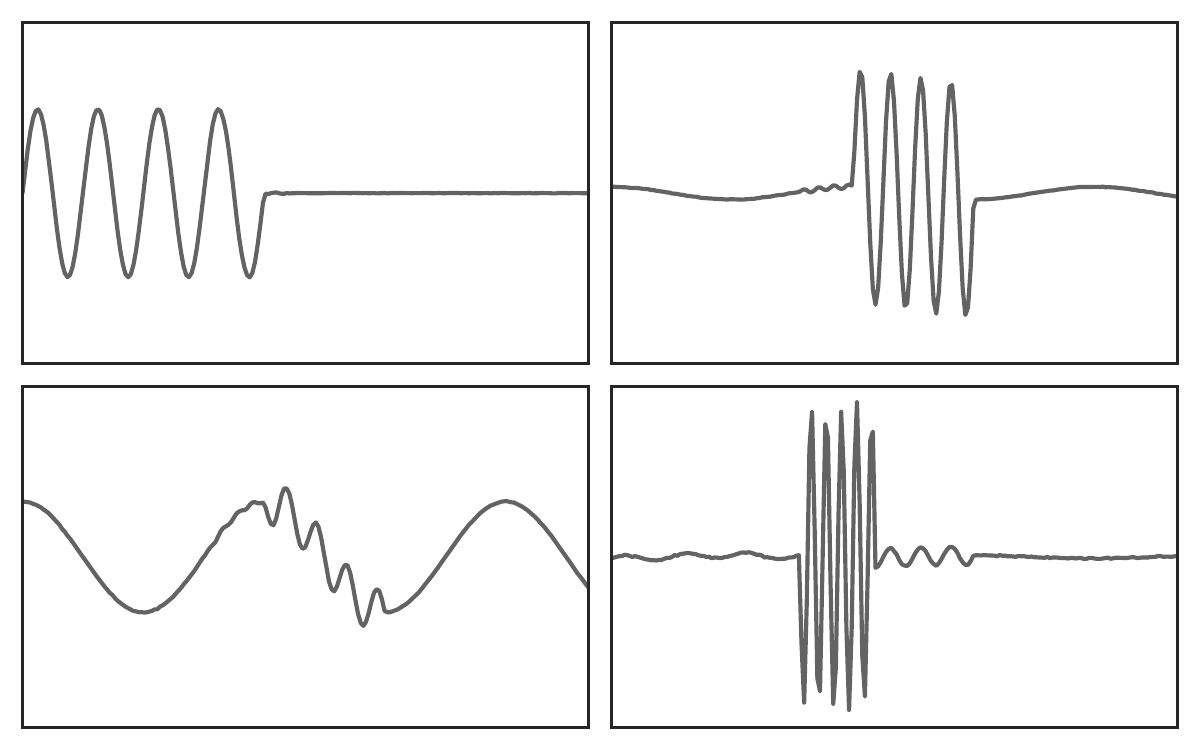}
		\caption{SVD}		
	\end{subfigure}
	
	\caption{Toy video decomposition results. The randomized CP with power iterations ($q=2$) accurately reconstructs the original spatiotemporal dynamics from noise-corrupted data, while SVD and rCP without power iterations ($q=0$) yield poor reconstruction results.}
	\label{fig:toy_video_decomposition}
\end{figure}

\begin{table}[!b]
	\centering
	\scalebox{0.85}{
		\begin{tabular}{ l l c c c c} 
			\hline 			\hline
			& \multicolumn{1}{l}{\bf Parameters}
			& \multicolumn{1}{c}{\bf Time (s)}
			& \multicolumn{1}{c}{\bf Speedup}
			& \multicolumn{1}{c}{\bf Iterations}
			& \multicolumn{1}{c}{\bf Error}									
			\\
			\cmidrule(r){1-6}
			
			\multirow{1}{*}{\rotatebox[origin=c]{0}{ \parbox{2.5cm}{CP BCD}  }} 
			& $k=4$ 		& 2.31	&  -  &  9  &  0.0171 \\ 
			\hline			
			
			\multirow{3}{*}{\rotatebox[origin=c]{0}{ \parbox{2.5cm}{rCP BCD} }} 
			& $k=4$, $p=10$, $q=0$              &  0.13	&  17   &  9  & 0.494   \\ 
			& $k=4$, $p=10$, $q=1$              &  0.14	&  16   & 10  &  0.0191   \\ 		
			& $k=4$, $p=10$, $q=2$              &  0.15	&  15   & 10  &  0.0164   \\ 
			\hline
			
			\multirow{1}{*}{\rotatebox[origin=c]{0}{ \parbox{2.5cm}{SVD}   }} 
			& $k=4$             &  0.52	&  -   &  -  &  0.137   \\
			\hline \hline
		\end{tabular}
	}
	\caption{Summary of the computational results for the noisy toy video.}
	\label{Tab:toy}
\end{table}

Figure~\ref{fig:toy_video_decomposition} displays the decomposition results of the noisy (SNR=$2$) toy video for both the randomized CP decomposition and the SVD. The first subplot shows the results of a rank $k=4$ approximation computed using the rCP algorithm with $q=2$ power iterations, and a small oversampling parameter $p=10$. The method faithfully captures the underlying spatial modes and the time dynamics. For illustration, the second subplot shows the decomposition results without additional power iterations. It can clearly be seen that this approach introduces distinct artifacts, and the approximation quality is relatively poor. 
The SVD, shown in the last subplot, demonstrates poor performance at separating the modes and mixes the spatiotemporal dynamics of modes 2 \& 3.

Table~\ref{Tab:toy} further quantifies the observed results. Interestingly, the relative error using the randomized algorithm with $q=2$ power iterations is slightly better than the deterministic algorithm. This is due to the intrinsic regularizing behavior of randomization. However, the reconstruction error without power iterations is large, as is the error resulting from the SVD. 
The achieved speedup of randomized CP is substantial, with a speedup factor of about $15$.
%
%
%

\paragraph{A Note on Compression.}
The CP decomposition provides a more parsimonious representation of the data.
Comparing the compression ratios between the CP decomposition and SVD illustrates the difference. For a rank $R=4$ tensor of dimension $100\times 100 \times 100$, the compression ratios are
\begin{eqnarray*}
	c_{SVD} & = &  \dfrac{I\cdot J\cdot K}{R\cdot (I\cdot J + K + 1)}  =  \dfrac{100^3}{4\cdot(100^2+ 100 + 1) } \approx 24.75, \\
	c_{CP} & = & \dfrac{I\cdot J\cdot K}{R\cdot (I + J + K + 1)}  =  \dfrac{100^3}{4\cdot(100 + 100 + 100 + 1) }  \approx  830.56.
\end{eqnarray*}

Note, the SVD requires the tensor to be reshaped in some direction. 
The comparison illustrates the striking difference between compression ratios. It is evident that the CP decomposition requires computing many fewer coefficients in order to approximate the tensor. Thus, less storage is required to approximate the data. 
%
While the CP decomposition yields a parsimonious approximation that is more robust to noise, the advantage of the SVD is that data can be approximated with an accuracy as low as machine precision.

\begin{figure}[!b]
	\centering
	\DeclareGraphicsExtensions{.png}
	\begin{overpic}[width=0.4\textwidth]{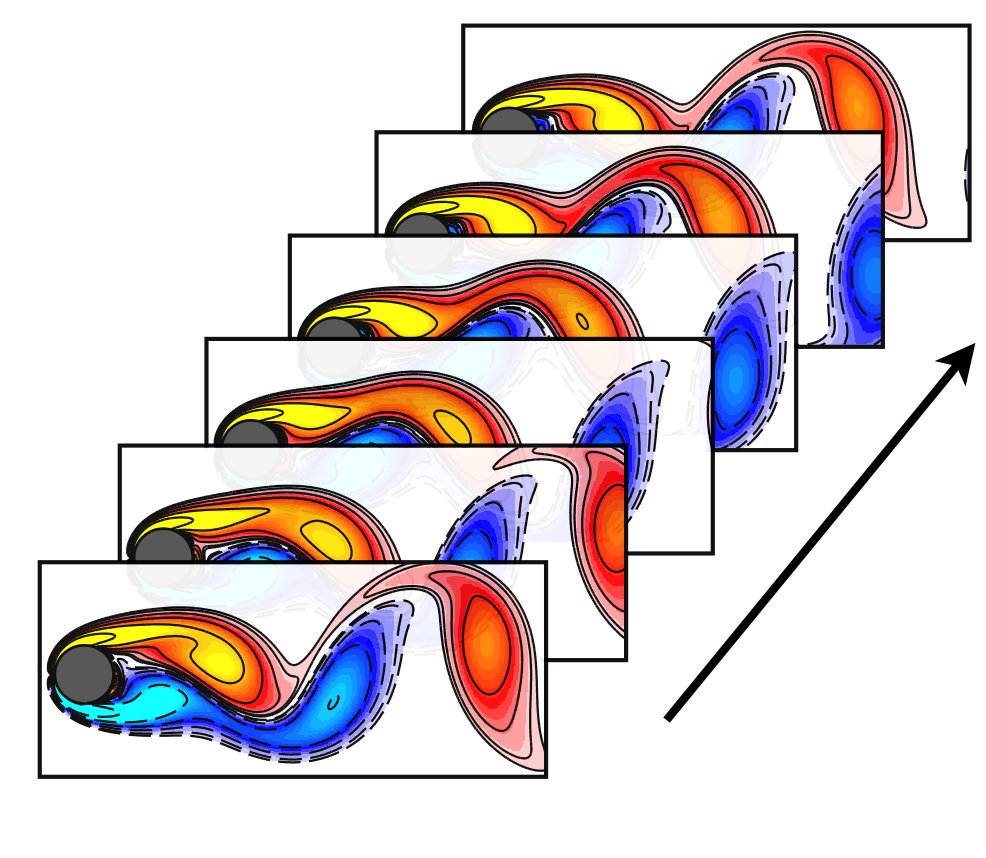} 
		\put(80,15){ time}
	\end{overpic}\vspace{-0.1cm}		
	
	\caption{Snapshots of the fluid flow behind a cylinder. }
	\label{fig:fluid_flow}
\end{figure}

\subsubsection{Flow behind a cylinder}

Extracting the dominant coherent structures from fluid flows helps to better characterize them for modeling and control~\citep{Brunton2015amr}. The workhorse algorithm in fluid dynamics and for model reduction is the SVD. However, fluid simulations generate high-resolution spatiotemporal grids of data which naturally manifest as tensors. In the following we examine the suitability of the CP decomposition for decomposing flow data, and compare the results to those of the SVD. 

The fluid flow behind a cylinder, a canonical example in fluid dynamics, is presented here. The data are a time-series generated from a simulation of fluid vorticity behind a stationary cylinder~\citep{colonius2008ibpm}. 
The corresponding flow tensor has dimension $199\times 449\times 151$, consisting of 151 snapshots on a $449\times 199$ spatial grid. Figure~\ref{fig:fluid_flow} shows a few example snapshots of the fluid flow.
%
%
The flow is characterized by a periodically shedding wake structure and is inherently low-rank in the absence of noise. 
%
%
%
The characteristic frequencies of flow oscillations occur in pairs, reflecting the complex-conjugate pairs of eigenvalues corresponding to harmonics in the temporal direction.
%
%

%
\begin{figure}[!t]
	\centering
	\begin{subfigure}[t]{0.41\textwidth}
		\centering
		\DeclareGraphicsExtensions{.png}
		\begin{overpic}[width=1\textwidth]{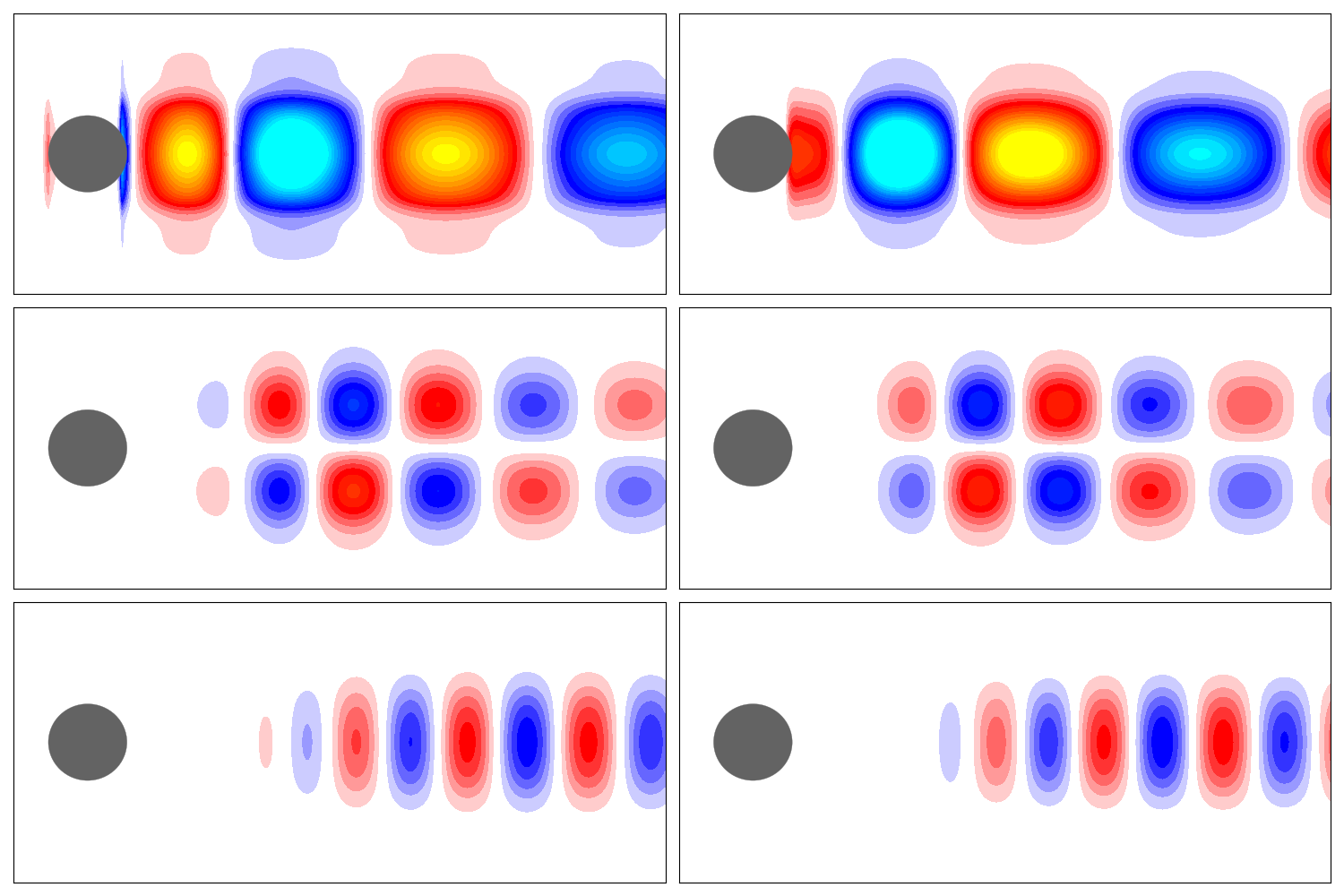}
			\put(-7,25){\rotatebox{90}{\small modes}}
		\end{overpic}		
		
	\end{subfigure}
	~
	\begin{subfigure}[t]{0.41\textwidth}
		\centering
		\DeclareGraphicsExtensions{.png}
		\includegraphics[width=1\textwidth]{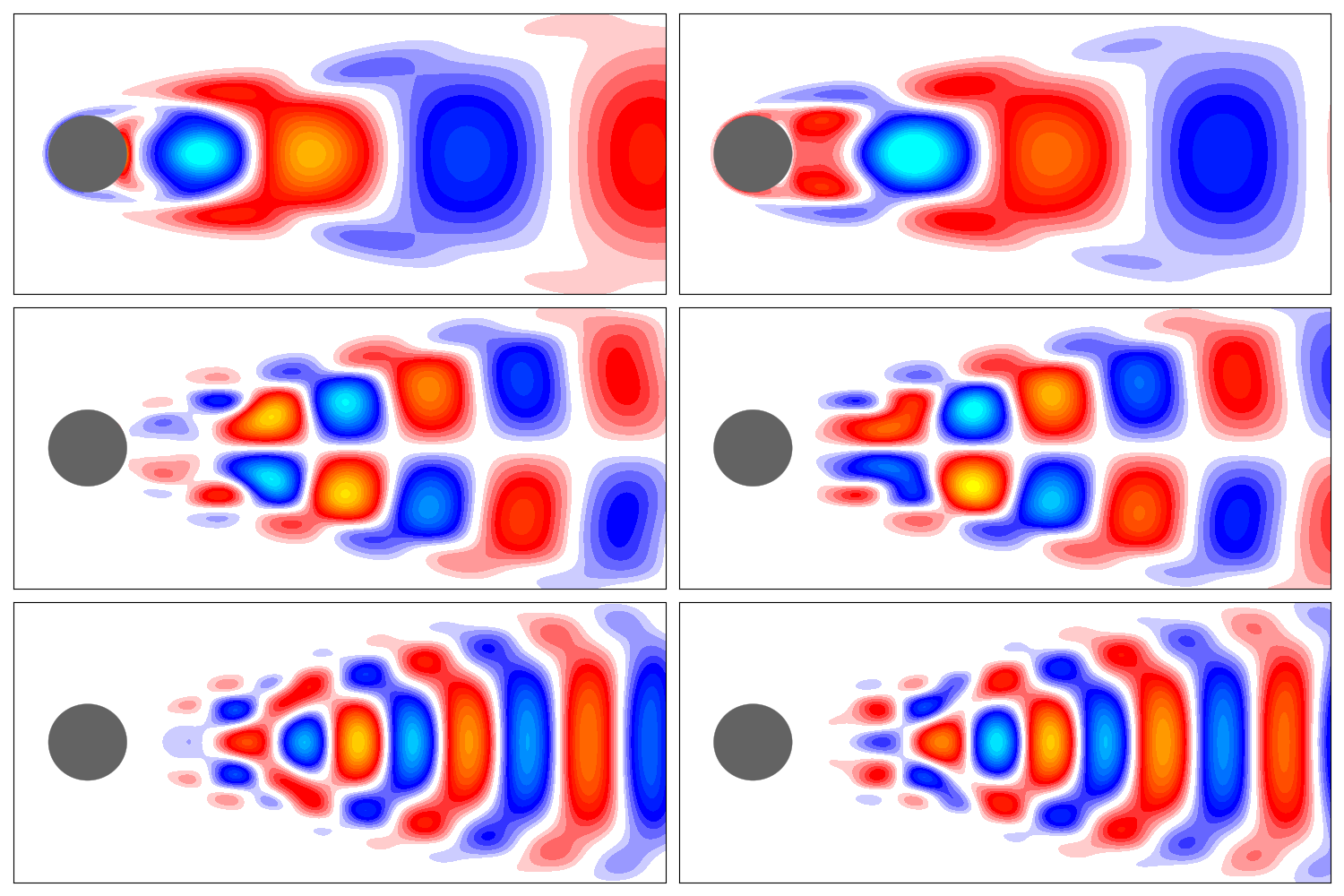}
		
	\end{subfigure}		
	
	\begin{subfigure}[t]{0.41\textwidth}
		\centering
		\DeclareGraphicsExtensions{.pdf}
		\begin{overpic}[width=1\textwidth]{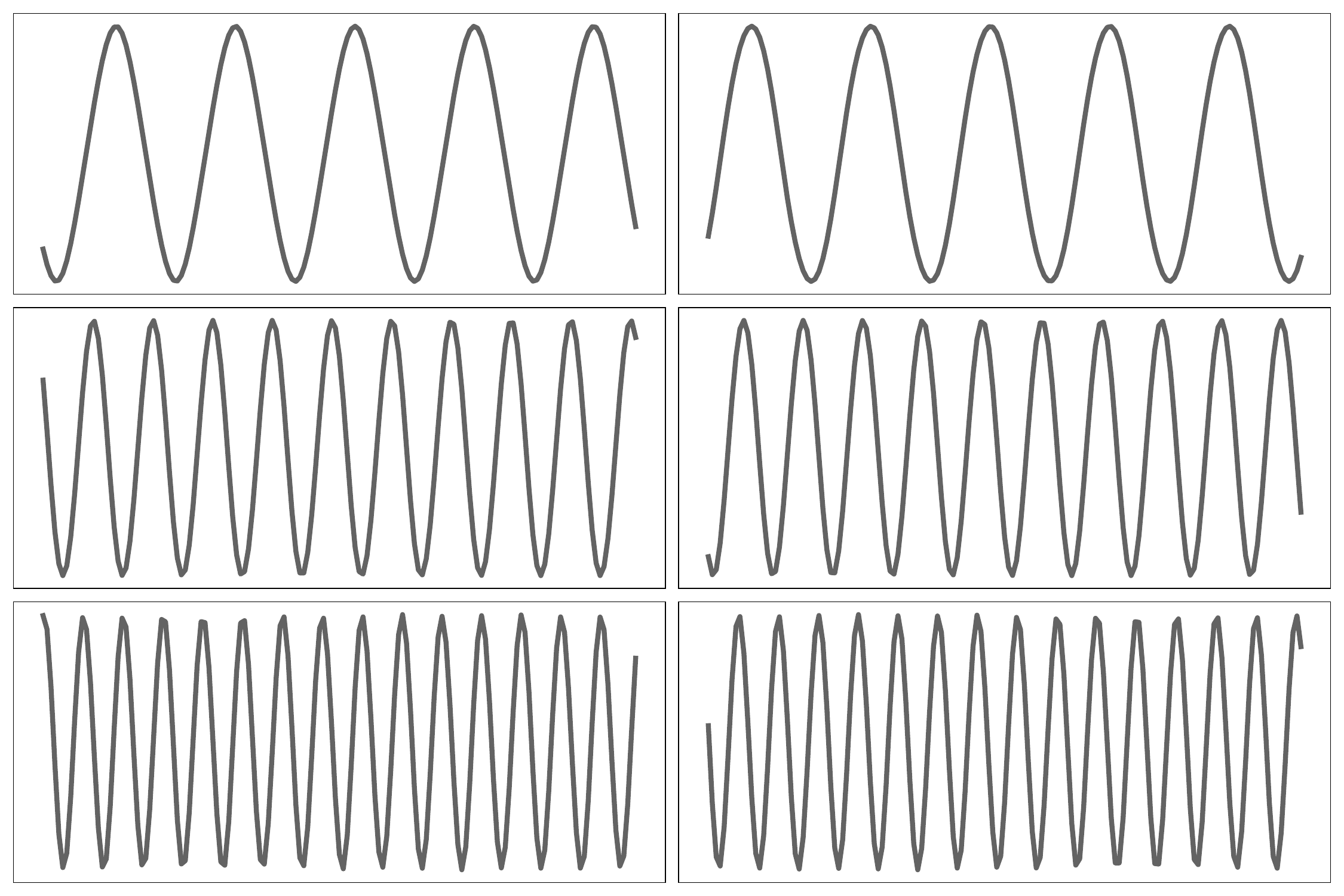}
			\put(-7,14){\rotatebox{90}{\small time dynamics}}
		\end{overpic}	
		\caption{Randomized CP ($q=2$).}		
	\end{subfigure}
	~
	\begin{subfigure}[t]{0.41\textwidth}
		\centering
		\DeclareGraphicsExtensions{.pdf}
		\includegraphics[width=1\textwidth]{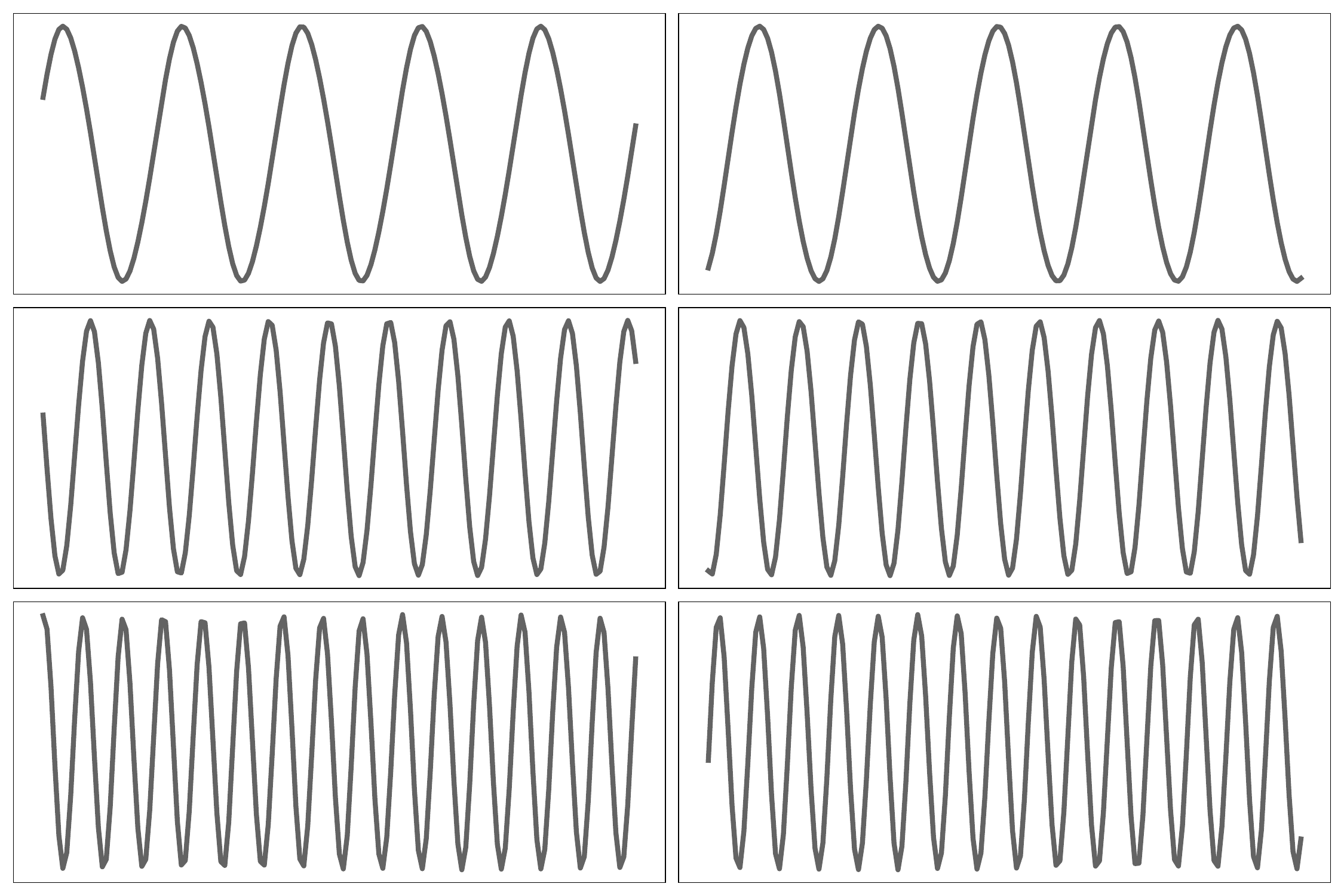}
		\caption{SVD.}
	\end{subfigure}	
	
	\caption{Fluid flow decomposition, no noise. Both methods capture the same dominant frequencies from the time dynamics, while randomized CP requires more rank-1 outer products in $x$ and $y$ to represent single frequency spatial dynamics.}
	\label{fig:fluid1}
\end{figure}

\begin{table}[!b]
	\centering
	\scalebox{0.85}{
		\begin{tabular}{ l l c c c c} 
			\hline 			\hline
			& \multicolumn{1}{l}{\bf Parameters}
			& \multicolumn{1}{c}{\bf Time (s)}
			& \multicolumn{1}{c}{\bf Speedup}
			& \multicolumn{1}{c}{\bf Iterations}
			& \multicolumn{1}{c}{\bf Error}									
			\\
			\cmidrule(r){1-6}
			
			\multirow{1}{*}{\rotatebox[origin=c]{0}{ \parbox{2.5cm}{CP BCD}  }} 
			& $k=30$ 		& 115.55	&  -  &  458  &  0.117 \\ 
			\hline			
			
			\multirow{3}{*}{\rotatebox[origin=c]{0}{ \parbox{2.5cm}{rCP BCD} }} 
			& $k=30$, $p=10$, $q=0$              &  1.27	&  91   &  533  & 0.122   \\ 
			& $k=30$, $p=10$, $q=1$              &  1.41	&  82   & 517  &  0.121   \\ 		
			& $k=30$, $p=10$, $q=2$              &  1.56	&  74   & 437  &  0.118   \\ 
			\hline
			
			\multirow{1}{*}{\rotatebox[origin=c]{0}{ \parbox{2.5cm}{SVD}   }} 
			& $k=30$             &  0.57	&  -   &  -  &  4.25E-05   \\
			\hline \hline
		\end{tabular}
	}
	\caption{Summary of the computational results for the noise-free cylinder flow.}
	\label{Tab:fluid}
\end{table}

\paragraph{Results in Absence of White Noise.}
Figure~\ref{fig:fluid1} shows both the approximated spatial modes and the temporal dynamics for the randomized CP decomposition and the SVD. 
Observe that both methods extract similar patterns, although unlike SVD, rCP spatial modes are sparse on the domain. The reason is two-fold: (i) CP does not impose orthogonality constraints (which in general reveal dense structure), and (ii) CP imposes rank-one outer product structure in the $x$ and $y$ directions via the columns of the factor matrices. In doing so, CP isolates the contributions of single wavenumbers (spatial frequencies) to the steady-state vortex shedding regime. These are commonly analyzed to study energy transfer between scales in complex flows and turbulence~\citep{meneveau2000scale,sharma2013coherent}, hence CP can help reveal new physically meaningful insights. In particular, randomization enables tensor decomposition of high-resolution vector fields that otherwise may not be tractable using deterministic methods. Thus rCP provides a novel decomposition of multimodal {\em fields} for downstream tasks in model reduction, pattern extraction and control.
%
%

Note, that for a fixed target rank of $k=30$ across all methods, the SVD achieves a substantially lower reconstruction error (see Table~\ref{Tab:fluid}). 
However, the compression ratios for the CP and SVD methods are $c_{CP}\approx 562.17$ and $c_{SVD} \approx 5.02$, i.e., the CP compresses the data nearly two orders of magnitude more. 

\begin{figure}[!t]
	\centering
	\begin{subfigure}[t]{0.41\textwidth}
		\centering
		\DeclareGraphicsExtensions{.png}
		\begin{overpic}[width=1\textwidth]{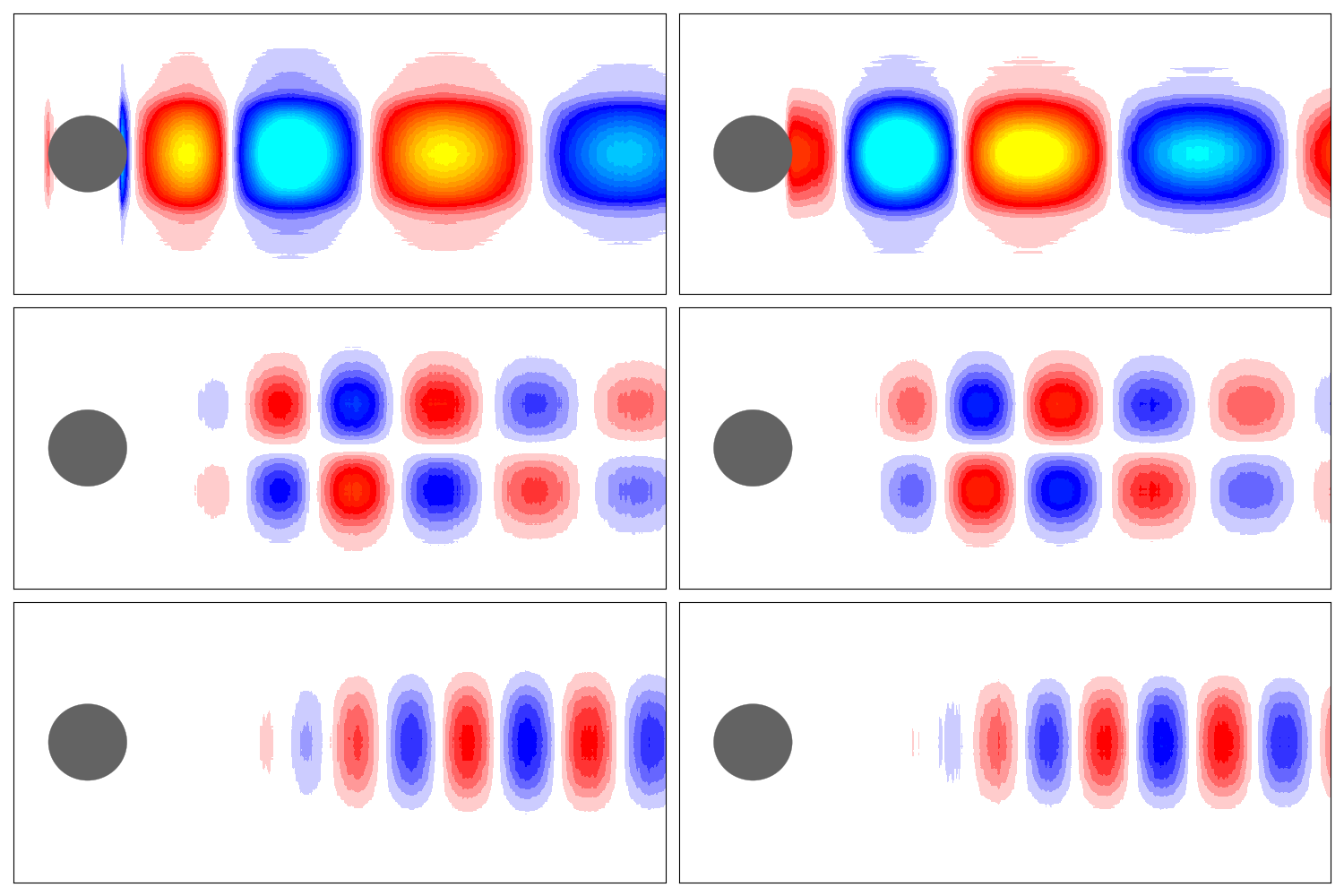}
			\put(-7,25){\rotatebox{90}{\small modes}}
		\end{overpic}			
		
	\end{subfigure}
	~
	\begin{subfigure}[t]{0.41\textwidth}
		\centering
		\DeclareGraphicsExtensions{.png}
		\includegraphics[width=1\textwidth]{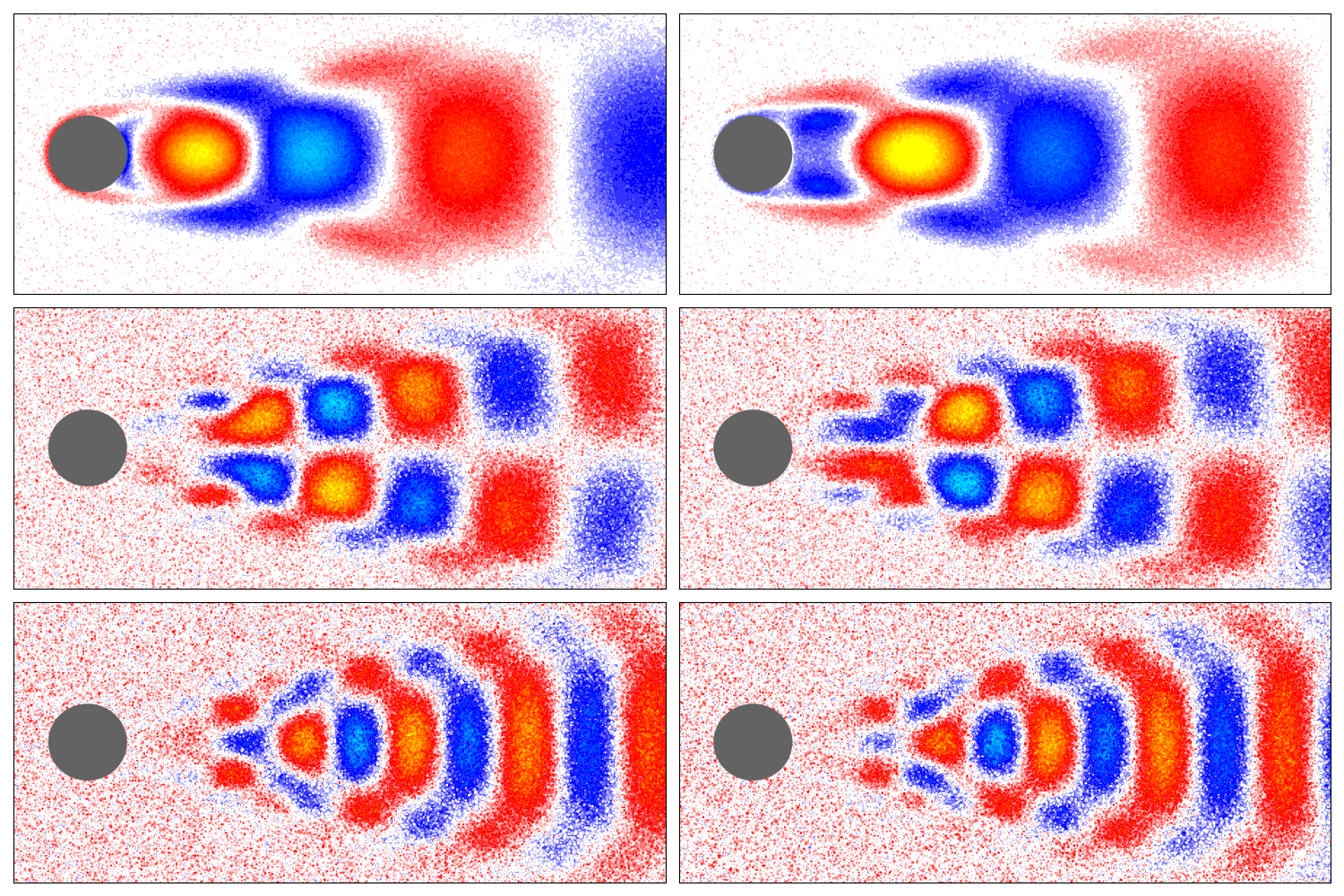}
	\end{subfigure}

	\begin{subfigure}[t]{0.41\textwidth}
		\centering
		\DeclareGraphicsExtensions{.pdf}
		\begin{overpic}[width=1\textwidth]{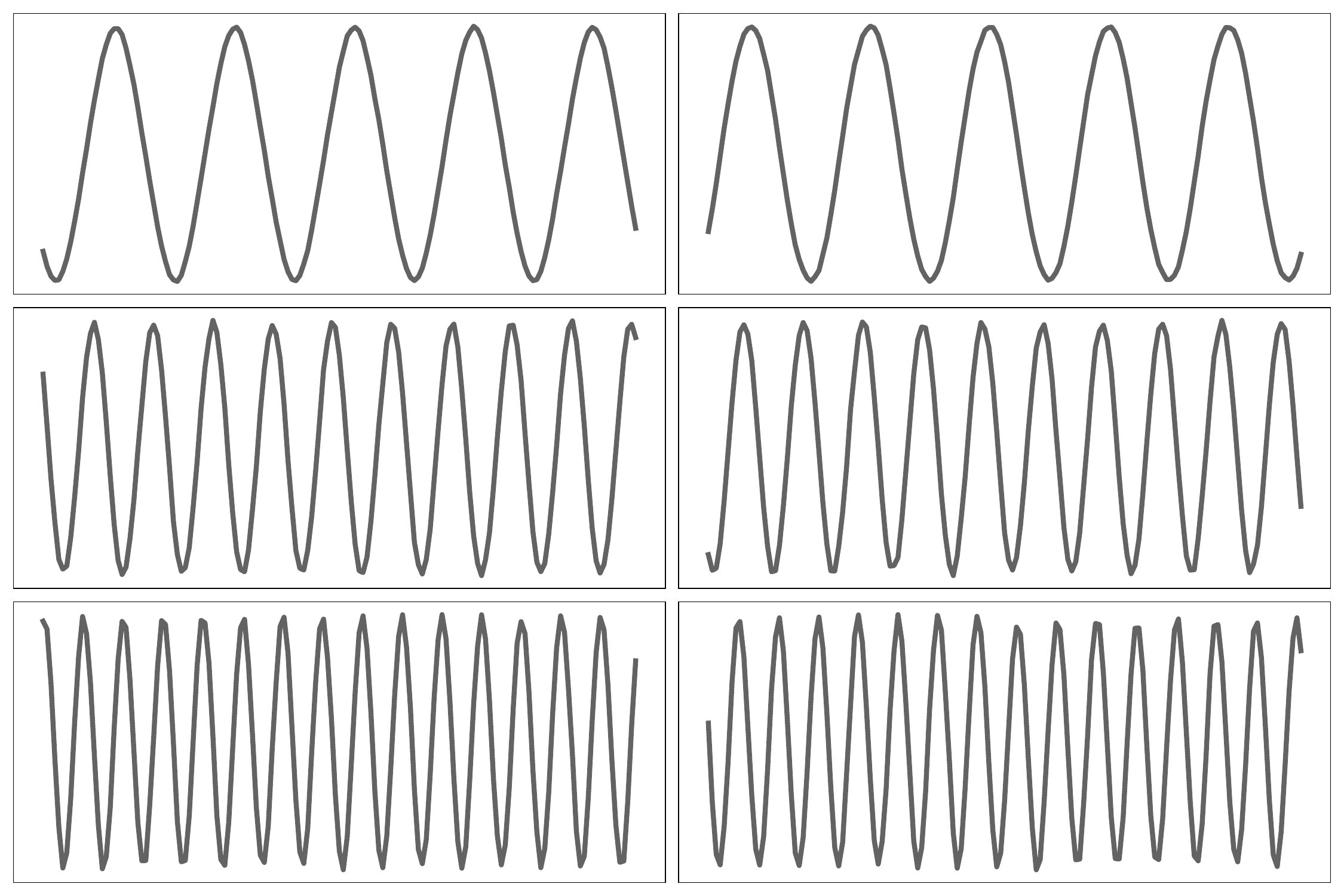}
			\put(-7,14){\rotatebox{90}{\small time dynamics}}
		\end{overpic}
		\caption{Randomized CP ($q=2$).}
	\end{subfigure}
	~
	\begin{subfigure}[t]{0.41\textwidth}
		\centering
		\DeclareGraphicsExtensions{.pdf}
		\includegraphics[width=1\textwidth]{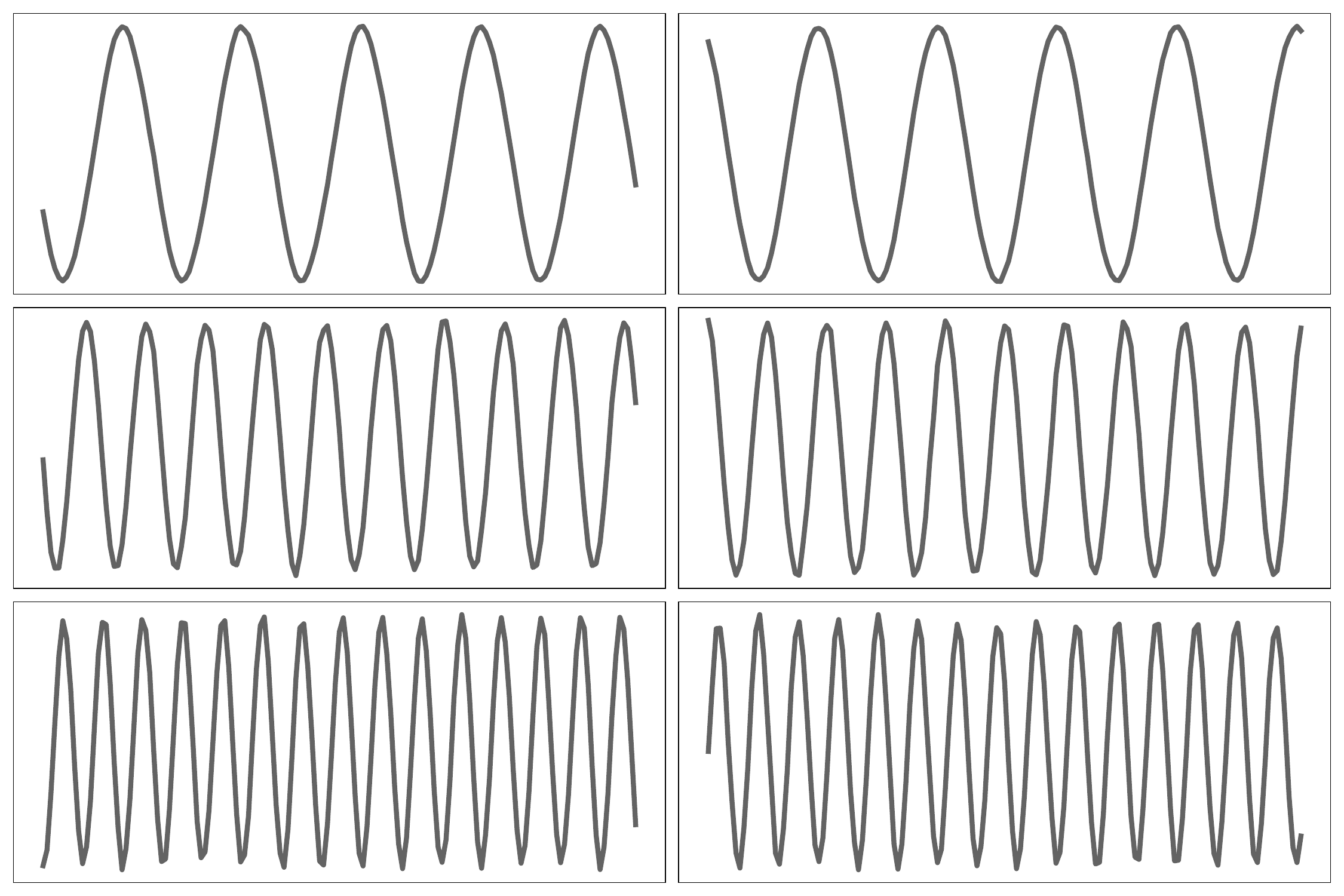}
		\caption{SVD.}
	\end{subfigure}	
	\caption{Fluid flow decomposition noisy (SNR=2). Randomized CP modes are robust to additive noise and SVD spatial modes are corrupted by noise.}
	\label{fig:fluid_noisy}
\end{figure}

\paragraph{Results in Presence of White Noise.}
Next, the analysis of the same flow is repeated in the presence of additive white noise. While this is not of concern when dealing with flow simulations, it is realistic when dealing with flows obtained from measurement. We chose a signal-to-noise ratio of 2 to demonstrate the robustness of the CP decomposition to noise.
Figure~\ref{fig:fluid_noisy} shows again the corresponding dominant spatial modes and temporal dynamics. Both the SVD and the CP decomposition faithfully capture the temporal dynamics. However, comparing the modes of the SVD to Figure~\ref{fig:fluid1}, it is apparent that the spatial modes contain a large amount of noise. The spatial modes revealed by the CP decomposition provide a significantly better approximation. Again, it is crucial to use power iterations to achieve a good approximation quality (see Table~\ref{Tab:fluid_noisy}). By inspection, the relative reconstruction error using the SVD is poor compared to the error achieved using the CP decomposition. Here, we show the error for a rank $k=30$ and $k=6$ approximation. The target rank was determined using the optimal hard threshold for singular values~\citep{gavish2014optimal}.

\begin{table}[!t]
	\centering
	\scalebox{0.85}{
		\begin{tabular}{ l l c c c c} 
			\hline 			\hline
			& \multicolumn{1}{l}{\bf Parameters}
			& \multicolumn{1}{c}{\bf Time (s)}
			& \multicolumn{1}{c}{\bf Speedup}
			& \multicolumn{1}{c}{\bf Iterations}
			& \multicolumn{1}{c}{\bf Error}									
			\\
			\cmidrule(r){1-6}
			
			\multirow{1}{*}{\rotatebox[origin=c]{0}{ \parbox{2.5cm}{CP BCD}  }} 
			& $k=30$ 		&  64.01	&  -  &  239  &  0.191 \\ 
			\hline			
			
			\multirow{3}{*}{\rotatebox[origin=c]{0}{ \parbox{2.5cm}{rCP BCD} }} 
			& $k=30$, $p=10$, $q=0$              &  0.99	&  64   &  332  & 0.522   \\ 
			& $k=30$, $p=10$, $q=1$              &  1.23	&  52   & 414  &  0.189   \\ 		
			& $k=30$, $p=10$, $q=2$              &  1.13	&  56   & 370  &  0.153   \\ 
			\hline
			
			\multirow{2}{*}{\rotatebox[origin=c]{0}{ \parbox{2.5cm}{SVD}   }} 
			& $k=30$             &  0.58	&  -   &  -  &  0.655   \\
			& $k=6$              &  0.58	&  -   &  -  &  0.311   \\							
			\hline \hline
		\end{tabular}
	}
	\caption{Summary of the computational results for the noise-corrupted cylinder flow.}
	\label{Tab:fluid_noisy}
\end{table}

The CP decomposition overcomes this disadvantage, and is able to approximate the first $k=30$ modes with only a slight loss of accuracy. Note that here the randomized CP decomposition performs better than the deterministic algorithm. We assume that this is due to the favorable intrinsic regularization effect of randomized methods.

\section{Conclusion} \label{sec:conclusion}

The emergence of massive tensors require efficient algorithms for obtaining tensor decompositions. 
To address this challenge, we have presented a randomized algorithm which substantially reduces the computational demands of the CP decomposition. 
Indeed, randomized algorithms have established themselves as highly competitive methods for computing traditional matrix decompositions. A key advantage of the randomized algorithm is that modern computational architectures are fully exploited. Thus, the algorithm benefits substantially from multithreading in a multi-core processor. In contrast to previously proposed high-performance tensor algorithms which are based on computational concepts such as distributed computing, our proposed randomized algorithm provides substantial computational speedups even on standard desktop computers. 
%
%
%
Our proposed algorithm achieves these speedups by reducing the computational costs per iteration, which enables the user to decompose real world examples that typically require a large number of iterations to converge.

In addition to computational savings, the randomized CP decomposition demonstrates outstanding performance on several examples using artificial and real-world data, including decompositions of high-resolution flow fields that may not be tractable with deterministic methods.
Moreover, our experiments show that the power iteration concept is crucial in order to achieve a robust tensor decomposition.  Thus, our algorithm has a practical advantage over previous randomized tensor algorithms, at a slightly increased computational cost due to additional power iterations.

\acks{We would like to thank Alex Williams, and Michael W. Mahoney for insightful discussion on tensor decompositions and randomized numerical linear algebra. Further, we would like to express our gratitude to the two anonymous reviewers for their valuable feedback, which helped us greatly improve the manuscript. NBE would like to acknowledge DARPA and NSF for providing partial support of this work. KM acknowledges support from NSF MSPRF Award 1803663. JNK acknowledges support from the Air Force Office of Scientific Research (AFOSR) grant FA9550-17-1-0329. SLB acknowledges funding support from the Air Force Office of Scientific Research (AFOSR) grant FA9550-18-1-0200.}

\section*{Data availability statement}

The data that support the findings of this study are available upon request.

\appendix
\section{Proof of Theorem 1} \label{app:proof1}
In the following, we sketch a proof for Theorem~\ref{thm:ubound}, which yields an upper bound for the approximate basis for the range of a tensor. 
To assess the quality of the basis matrices $\{\mathbf{Q}_n\}_{n=1}^{N}$, we first show that the problem can be expressed as a sum of subproblems. Defining the residual error 
\begin{equation}\label{eq:residual}
\|\tE \|_F = \|  \tX - \hat{ \tX } \| =  \| \tX - \tX \times_1 \bQ_1\bQ_1^\top \times_2 \cdots \times_N \bQ_N\bQ_N^\top \|_F.
\end{equation}
Note that the Frobenius norm of a tensor and its matricized forms are equivalent. Defining the orthogonal projector $\mathbf{P}_n \equiv \bQ_n\bQ_n^\top $, we can reformulate \eqref{eq:residual}  compactly as
\begin{equation}\label{eq:residual2}
\|\tE \|_F = \|\tX - \tX \times_1  \mathbf{P}_1 \times_2 \cdots \times_N \mathbf{P}_N    \|_F.
\end{equation}
\begin{proof} Assuming that $\mathbf{P}_n$ yields an exact projection onto the column space of the matrix $\mathbf{Q}_n$, we need to show first that the error can be expressed as a sum of the errors of the $n$ projections
	\begin{equation}\label{eq:residualSum}
	\|\tE \|_F = \sum_{n=1}^{N} \|  \tX  - \tX \times_n \mathbf{P}_n \|_F = \sum_{n=1}^{N} \|  \tX \times_n (\mathbf{I} - \mathbf{P}_n) \|_F,
	\end{equation}
	where $\mathbf{I}$ denotes the identity matrix. Following~\cite{DrineasTensor}, let us add and subtract the term $\tX \times_N \mathbf{P}_N$ in Equation~\eqref{eq:residual2} so that we obtain
	\begin{eqnarray}
	\|\tE \|_F & = & \|\tX - \tX \times_N \mathbf{P}_N + \tX \times_N \mathbf{P}_N - \tX  \times_1 \mathbf{P}_1 \times_1 \cdots \times_N \mathbf{P}_N \|_F \\
	& \leq & \|\tX - \tX \times_N \mathbf{P}_N \|_F + \|\tX \times_N \mathbf{P}_N - \tX  \times_1 \mathbf{P}_1 \times_1 \cdots \times_N \mathbf{P}_N  \|_F \label{eq:tri2}\\ 
	& = & \|\tX - \tX \times_N \mathbf{P}_N \|_F + \|(\tX - \tX  \times_1 \mathbf{P}_1 \times_1 \cdots \times_{N-1} \mathbf{P}_{N-1}) \times_N \mathbf{P}_N  \|_F \label{eq:tri3}\\  
	& \leq & \|\tX - \tX \times_N \mathbf{P}_N \|_F + \|\tX - \tX  \times_1 \mathbf{P}_1 \times_1 \cdots \times_{N-1} \mathbf{P}_{N-1}  \|_F. \label{eq:tri4}
	\end{eqnarray}
	The bound~\eqref{eq:tri2} follows from the triangular inequality for a norm. Next, the common term $\mathbf{P}_N$ is factored out in Equation~\eqref{eq:tri3}. Then, the bound~\eqref{eq:tri4} follows from the properties of orthogonal projectors. This is because the $range(\tX  \times_1 \mathbf{P}_1 \times_1 \cdots \times_{N-1}  \mathbf{P}_{N-1}) \subset range(\tX  \times_1 \mathbf{P}_1 \times_1 \cdots \times_N \mathbf{P}_N)$, and then it holds that $\|\tX - \tX  \times_1 \mathbf{P}_1 \times_1 \cdots \times_N \mathbf{P}_N \|_F \leq \|\tX - \tX  \times_1 \mathbf{P}_1 \times_1 \cdots \times_{N-1}  \mathbf{P}_{N-1} \|_F$. See Proposition 8.5 by \cite{halko2011rand} for a proof using matrices. Subsequently the residual error  $\tE_{N-1}$ can be bounded
	\begin{eqnarray}
	\|\tE_{N-1}\| & \leq & \| \tX - \tX \times_{N-1} \mathbf{P}_{N-1}   \|_F + \| \tX - \tX  \times_1 \mathbf{P}_1 \times_1 \cdots \times_{N-2} \mathbf{P}_{N-2})  \|_F.
	\end{eqnarray}
	From this inequality, Equation~\eqref{eq:residualSum} follows.
	We take the expectation of Equation~\eqref{eq:residualSum} 
	\begin{equation}\label{eq:EresidualSum}
	\E\|\tE \|_F  = \E \left[ \sum_{n=1}^{N} \|  \tX \times_n (\mathbf{I} - \mathbf{P}_n) \|_F \right].
	\end{equation}
	Recalling that Theorem 10.5 formulated by~\cite{halko2011rand} states the following expected approximation error (formulated here using tensor notation)
	\begin{eqnarray}\label{eq:halkoThem}
	\E\| \tX \times_n (\mathbf{I} - \mathbf{P}_n) \|_F \leq \sqrt{ 1 + \frac{k}{p-1}} \cdot \sqrt{\sum_{j>k}\sigma_{j}^2},
	\end{eqnarray}
	assuming that the sketch in Equation~\eqref{eq:sampleMatrix} is constructed using a standard Gaussian matrix $\mathbf{\Omega}$. Here $\sigma_{j}$ denotes the singular values of the matricized tensor $\tX_{(n)}$ greater then the chosen target rank $k$. Combining Equations~\eqref{eq:EresidualSum} and~\eqref{eq:halkoThem} then yields the results of the theorem (\ref{eq:theorem}).
\end{proof} 

\begin{figure}[H]
	\centering
	\begin{subfigure}[t]{0.75\textwidth}
		\centering
		\DeclareGraphicsExtensions{.png}
		\includegraphics[width=1\textwidth]{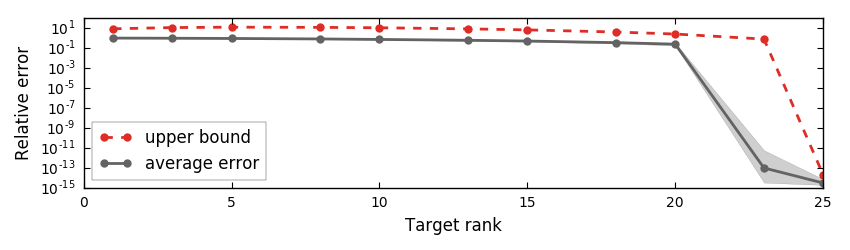}
		\caption{Tensor of dimension $50\times 50\times 50$.}
	\end{subfigure}
	
	\begin{subfigure}[t]{0.75\textwidth}
		\centering
		\DeclareGraphicsExtensions{.png}
		\includegraphics[width=1\textwidth]{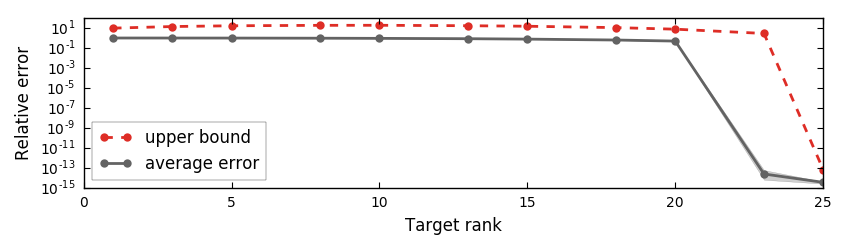}
		\caption{Tensor of dimension $50\times 50\times 50\times 50$.}
	\end{subfigure}
	%
	\caption{Empirical evaluation of the theoretical upper bound. }
	\label{fig:validation}
\end{figure}

Figures~\ref{fig:validation} evaluates the theoretical upper bound over 100 runs for both a third and fourth order random low-rank $R=25$ tensor. Here, we use a fixed oversampling parameter $p=2$. The results show that the empirical error is faithfully bounded by the theoretical upper bound for varying target ranks.

\small
\bibliography{dis}

\end{document}